\documentclass{amsart}
\usepackage{amsmath,amssymb}

\newtheorem{theorem}{Theorem}[section]                    
\newtheorem{proposition}[theorem]{Proposition}            
\newtheorem{corollary}[theorem]{Corollary}                
\newtheorem{lemma}[theorem]{Lemma}
\newtheorem{remark}[theorem]{Remark} 
\newtheorem{remarks}[theorem]{Remarks}

\newtheorem{fact}[theorem]{Fact}
\newtheorem{facts}[theorem]{Facts}
\newtheorem{definition}{Definition}[section]

\begin{document}

\title[HNN Extensions]
{HNN Extensions of von Neumann Algebras}
\author[Y. Ueda] {Yoshimichi UEDA} 
\address{Graduate School of Mathematics, Kyushu University, Fukuoka, 810-8560, Japan}
\email{ueda@math.kyushu-u.ac.jp}
\thanks{Supported by Grant-in-Aid for Young Scientists (B) 14740118.}
\maketitle

\begin{abstract} 
Reduced HNN extensions of von Neumann algebras (as well as $C^*$-algebras) will be introduced, and their modular theory, factoriality and ultraproducts will be discussed. 
In several concrete settings, detailed analysis on them will be also carried out.    
\end{abstract}   

\section{Introduction} 
There are two fundamental constructions in combinatorial or geometric group theory, which are those of free products with amalgamations and of HNN (G.~Higman, B.~H.~Neumann and H.~Neumann \cite{higmanneumann^2:JLondonMath1949}) extensions. The interested reader may consult \cite{lyndonschupp:book} as a standard reference on the topics. Even in the framework of von Neumann algebras (as well as $C^*$-algebras), reduced free products with amalgamations  (\cite{voiculescu:lecturenotes1132}\cite{voiculescudykemanica:monograph}\cite{popa:invent:irredindex>4} and also \cite{ueda:pacific}) have been seriously investigated so far and played key r\^{o}les in several resolutions of ``existence" questions in the theory of von Neumann algebras (see, e.g. \cite{popa:invent:irredindex>4}\cite{radulescu:freegroupsubfactor}\cite{shlyakhtenkoueda:crelle}\cite{popashlyakhtenko:preprint} and also \cite{ueda:journalofmathsocJpn}). However, HNN extensions have never been discussed so far in the framework. 

Historically, many ideas in group theory, especially part of dealing with countably infinite discrete groups, have been applied directly and/or indirectly to many aspects in the theory of von Neumann algebras (as well as $C^*$-algebras) since the beginning of the theory. In fact, many explicit examples of von Neumann algebras that opened new perspectives in the theory came from group theory (see e.g. \cite{murrayvonneumann4}\cite{mcduff:uncountableexamples}\cite{haagerup:invent79}\cite{connes:jot80}\cite{connes-jones:bull-london}\cite{cowling-haagerup:invent89} and also recent breakthroughs \cite{popa:betti}\cite{ozawa:procAMS}\cite{ozawa:acta}\cite{ozawapopa:invent}), and it is still expected to find much more ``monsters" (i.e., concrete examples with very special properties) living in the world of non-amenable von Neumann algebras. To do so, it seems still to be one of the important guiding principles to seek for new ideas in group theory. Following this principle, we will introduce reduced HNN extensions in the framework of von Neumann algebras (as well as $C^*$-algebras) and take a very first step towards serious and systematic investigation on them with aiming that their construction will play a key r\^{o}le in future attempts of constructing new monsters in the world of non-amenable von Neumann algebras. 

Let us explain the organization of this article. In \S2, we will review free products with amalgamations of von Neumann algebras with special emphasis of the admissibility of embedding maps of amalgamated algebras in the construction. Although this slight generalization of the previously used one is of course a folklore, we will briefly review it to avoid any confusion since the admissibility of embedding maps plays a key r\^{o}le in our construction of HNN extensions. In \S3, reduced HNN extensions of von Neumann algebras will be introduced, and then their characterization (or their ``construction-free" definition) given in terms of expected algebraic relations and ``moment-values" of conditional expectations as in the case of free products with amalgamations. In the group setting, one standard way of constructing HNN extensions is the use of ``shift automorphisms" on ``infinite free products with amalgamations" over {\it isomorphic but not necessary common} subgroups (in fact, two different embeddings of amalgamated groups are needed). This amalgamation procedure brings us ``difficulty" in constructing ``shift automorphisms" in connection with conditional expectations since the universal construction is not applicable in the von Neumann algebra setting. Hence, a different idea is needed to construct the desired ones, and indeed it is based on an observation coming from our previous work \cite{ueda:us-jpn} on a different topic. Roughly speaking, our construction can be understood as an ``amalgam" (but not a ``combination") of those of covariant representations without unitary implementations in the crossed-product construction (see \cite[Vol.II; Eq.~(10) in p.241]{takesaki:book}) and of free products with amalgamations. Our construction seems somewhat natural from the group theoretic viewpoint. In fact, the notion of HNN extensions is known to be necessary to describe a subgroup of a given free product group with amalgamation over a non-trivial subgroup. The \S4 will concern modular theoretical aspects of reduced HNN extensions. More precisely, we will give a complete description of modular automorphisms and also show that the continuous core of any reduced HNN extension becomes again a reduced HNN extension. In \S5, we will discuss the factoriality and investigate the ultraproducts of reduced HNN extensions. The results correspond to what we obtained in our previous work \cite{ueda:transAMS} on free products with amalgamations. In \S6, we will investigate reduced HNN extensions of von Neumann algebras in several concrete settings. The first one is naturally arisen from non-commutative $2$-tori, the second from the tensor product operation, and the third from regular and singular MASAs in the crossed-products by (non-commutative) Bernoulli shifts. The third one seems important for further investigation because any given surjective (partial) $*$-isomorphism between regular and singular MASAs in question can never be extended to any global $*$-automorphism on the given ``base" algebras. In \S7, reduced HNN extensions of $C^*$-algebras will be introduced in the same manner as in the von Neumann algebra setting, and then some basic facts will be given. Further analysis on them will be presented elsewhere. 
 

Part of this article was presented in the conference ``Recent Advances in von Neumann Algebras" celebrated to Professor Masamichi Takesaki's 70th birthday, at UCLA in May, 2003. We would like to express our sincere thanks to the organizers; Professors Yasuyuki Kawahigashi, Sorin Popa, and Dimitri Shlyakhtenko, who kindly gave us the opportunity to present this work in the conference, and also would like to celebrate Professor Masamichi Takesaki's 70th birthday.   

\medskip\noindent
{\it Acknowledgment.} We would like to express our  sincere appreciation to Yasuo Watatani who first reminded us, through many excursive conversations, HNN extensions of groups and the importance to look for ideas in geometric group theory, when we had worked out free products with amalgamations over Cartan subalgebras. We thank Tomohiro Hayashi for fruitful conversations at the final stage of this work, and also Masaki Izumi for his useful comment, to which the present form of Corollary \ref{Cor3-6} is indebted. We thank also the anonymous referee for valuable comments.   

\section{Preliminaries on Free Products with Amalgamations} 

Let $D$ and $N_s$ ($s \in S$, an index set) be $\sigma$-finite von Neumann algebras, and we have a normal $*$-isomorphism $\iota_s : D \rightarrow N_s$ for each $s \in S$. Suppose further that the von Neumann subalgebra $\iota_s(D)$ of $N_s$ is the range of a faithful normal conditional expectation $E_s$ for every $s \in S$. Even in this setting, we will be still able to construct the reduced free product with amalgamation 
\begin{equation*}
\left(N, E\right) = \underset{s \in S}{\bigstar_D} \left(N_s, E_s : \iota_s\right).  
\end{equation*}      
The discussions in this article will treat the type II and type III cases in common so that the approach in \cite{ueda:pacific} to the amalgamated free product construction will be convenient since complete treatment of modular theory was given there. To construct reduced HNN extensions, the admissibility of the embeddings $\iota_s$ ($s \in S$) in the construction plays a key r\^{o}le. Hence, following \cite{ueda:pacific} we would like to recall (without details) the amalgamated free product construction with special emphasis on the embeddings $\iota_s$'s to avoid any confusion. 

Fix $s \in S$ for a while, and let $\left({\mathcal H}_s, N_s, J_s, {\mathcal P}_s^{\natural}\right)$ and $\left(L^2(D), D, J_D, {\mathcal P}_D^{\natural}\right)$ be the standard forms. See \cite[Vol.II; Chap.~IX, \S1]{takesaki:book} for detailed account of standard forms. Using the mapping 
\begin{equation}\label{Eq-1}
\xi \in {\mathcal P}_D^{\natural} \longmapsto 
\left(\left(\omega_{\xi}\big|_D\right)\circ\iota_s^{-1}\circ E_s\right)^{\frac{1}{2}} 
\in {\mathcal P}_s^{\natural} 
\end{equation}
we can extend the embedding $\iota_s : D \rightarrow N_s$ to the Hilbert space level and still denote it by the same symbol $\iota_s : L^2(D) \rightarrow {\mathcal H}_s$. Here, $\psi^{\frac{1}{2}} \in {\mathcal P}_s^{\natural}$ denotes the unique implementing vector of a normal positive linear functional $\psi$ on $N_s$. This embedding satisfies the following expected properties: (i) For $\xi \in L^2(D)$ and $d_1, d_2 \in D$, we have $\iota_s\left(d_1 J_D d_2^* J_D \xi\right) = \iota_s\left(d_1\right) J_s\iota_s\left(d_2\right)^*J_s\iota_s\left(\xi\right)$, i.e., $\iota_s\left(d_1\cdot\xi\cdot d_2\right) = 
\iota_s\left(d_1\right)\cdot\iota_s\left(\xi\right)\cdot\iota_s\left(d_2\right)$ with the usual notations in the bimodule theory. (ii) For each $\xi \in {\mathcal P}_D^{\natural}$, the vector $\iota_s\left(\xi\right)$ becomes the canonical implementing one in ${\mathcal P}^{\natural}_s$ of the state $\left(\omega_{\xi}\big|_D\right)\circ\iota_s^{-1}\circ E_s$, a consequence from \eqref{Eq-1}.

Fix a faithful normal state $\varphi$ on $D$ and denote by $\xi_{\varphi}$ its implementing vector in ${\mathcal P}_D^{\natural}$. As mentioned above, the vector $\iota_s\left(\xi_{\varphi}\right)$ becomes the unique implementing one of the state $\varphi\circ\iota_s^{-1}\circ E_s$ in the natural cone ${\mathcal P}_s^{\natural}$. We denote the kernel of $E_s$ by $N_s^{\circ}$ as usual, and introduce the operation $x \in N_s \mapsto x^{\circ} := x - E_s(x) \in N_s^{\circ}$. We also write ${\mathcal H}_s^{\circ} := {\mathcal H}_s \ominus \iota_s\left(L^2(D)\right)$, and it is clear that the subspace ${\mathcal H}_s^{\circ}$ is invariant under the left and right actions of $D$ via the embedding map $\iota_s$.  Thus the natural $D$-$D$ bimodule structure of the Hilbert space ${\mathcal H_s}$: 
\begin{equation*}
d_1\cdot\xi\cdot d_2 := \iota_s\left(d_1\right) J_s\iota_s\left(d_2\right)^*J_s\xi, \quad \xi \in 
{\mathcal H}_s,\ d_1,d_2 \in D
\end{equation*}
is inherited to the subspace ${\mathcal H}_s^{\circ}$. When emphasize this bimodule structure, we will use the symbols ${}_D\left({}_{\iota_s}{\mathcal H}_s{}_{\iota_s}\right)_D$, ${}_D\left({}_{\iota_s}{\mathcal H}_s^{\circ}{}_{\iota_s}\right)_D$ (or ${}_{\iota_s}{\mathcal H}_s^{\circ}{}_{\iota_s}$, ${}_{\iota_s}{\mathcal H}_s^{\circ}{}_{\iota_s}$ for short) instead of ${\mathcal H}_s$, ${\mathcal H}_s^{\circ}$, respectively. Notice here that we have the natural bimodule isomorphism 
\begin{equation*}
{}_DL^2(D)_D\oplus{}_D\left({}_{\iota_s}{\mathcal H}_s^{\circ}{}_{\iota_s}\right)_D \cong {}_D\left({}_{\iota_s}{\mathcal H}_s{}_{\iota_s}\right)_D 
\end{equation*}
given by $\xi\oplus\eta \longmapsto \iota_s\left(\xi\right) + \eta$. Let us construct the Hilbert space 
\begin{equation*}
{\mathcal H} := 
L^2(D)\oplus\sideset{}{^{\oplus}}\sum_{s_1 \neq s_2 \neq \cdots \neq s_n} 
{}_{\iota_{s_1}}{\mathcal H}^{\circ}{}_{\iota_{s_1}}\otimes_{\varphi}
{}_{\iota_{s_2}}{\mathcal H}^{\circ}{}_{\iota_{s_2}}\otimes_{\varphi} \cdots \otimes_{\varphi} 
{}_{\iota_{s_n}}{\mathcal H}^{\circ}{}_{\iota_{s_n}}, 
\end{equation*}
on which the desired algebra $N$ acts. This naturally becomes a $D$-$D$ bimodule, and the left and right actions are denoted by $\lambda$ and $\rho$, respectively. For each $s \in S$, we can construct the $*$-representation $\lambda_s : N_s \rightarrow \mathrm{End}\left({\mathcal H}_D\right)$ and the anti-$*$-representation $\rho_s : N_s \rightarrow \mathrm{End}\left({}_D{\mathcal H}\right)$ by the same way as in \cite[p.361--362]{ueda:pacific}. To do so, we need only some basic properties on relative tensor products (see \cite[Vol.II; Chap.~IX, \S3]{takesaki:book}) and the bimodule isomorphism ${}_{\iota_s}{\mathcal H}_s{}_{\iota_s} \cong L^2(D)\oplus\left({}_{\iota_s}{\mathcal H}_s^{\circ}{}_{\iota_s}\right)$ precisely explained above. Let us consider two von Neumann algebras 
\begin{equation*}
N := \left(\bigcup_{s \in S} \lambda_s\left(N_s\right)\right)'', \quad 
L := \left(\bigcup_{s \in S} \rho_s\left(N_s\right)\right)'' \quad 
\text{on ${\mathcal H}$}, 
\end{equation*}
and define $\psi := \omega_{\xi_{\varphi}}\big|_N$, as a vector state, with regarding $\xi_{\varphi} \in L^2(D)$ as a vector in ${\mathcal H}$. 
      
\begin{facts} {\rm (\cite[p.362--365]{ueda:pacific})}\label{Facts2-1}
\begin{itemize}
\item[(A)] $\lambda_s\circ\iota_s$ coincides with the left action $\lambda$ of $D$ for each $s \in S$. 
\item[(B)] $\rho_s\circ\iota_s$ coincides with the right action $\rho$ of $D$ for each $s \in S$. 
\item[(C)] The vector $\xi_{\varphi}$ is cyclic for both $N$ and $L$. 

\item[(D)] The commutant $N'$ on ${\mathcal H}$ contains $L$. {\rm (}More on this is true, that is, the commutation theorem $N' = L$ holds, see \cite[Appendix II]{ueda:pacific}.{\rm )} Hence, the state $\psi$ is faithful. 
\item[(E)] For each $x_j^{\circ} \in N_{s_j}^{\circ}$ with $s_1 \neq s_2 \neq \cdots \neq s_n$, we have 
\begin{equation*}
\psi\left(\lambda_{s_1}\left(x_1^{\circ}\right)\cdots\lambda_{s_n}\left(x_n^{\circ}\right)\right) = 0. 
\end{equation*}
\item[(F)] The modular automorphism $\sigma_t^{\psi}$ {\rm (}$t \in {\mathbf R}${\rm )} 
satisfies 
\begin{equation*}
\sigma_t^{\psi} \circ \lambda_s = \lambda_s\circ\sigma_t^{\varphi\circ\iota_s^{-1}\circ E_s}\quad 
\text{and}\quad \sigma_t^{\psi}\circ\lambda = \lambda\circ\sigma_t^{\varphi}. 
\end{equation*}
Hence, there is a {\rm (}unique{\rm )} $\psi$-preserving conditional expectation $E^{\psi} : N \rightarrow \lambda(D)$ thanks to Takesaki's theorem {\rm (}\cite[Vol.II; Theorem 4.2 in Chap.~IX]{takesaki:book}{\rm )}. 
\end{itemize}
\end{facts}  

As in \cite[lines 8--3 from the bottom in p.364]{ueda:pacific}, the above (C),(E) imply the freeness (with amalgamation over $\lambda(D)$) among the von Neumann subalgebras $\lambda_s\left(N_s\right)$ ($s \in S$) with respect to $E^{\psi}$ in the sense of Voiculescu \cite[\S5]{voiculescu:lecturenotes1132}:
\begin{equation*}
E^{\psi}\left(\lambda_{s_1}\left(x_1^{\circ}\right)\cdots\lambda_{s_n}\left(x_n^{\circ}\right)\right) = 0
\end{equation*}
whenever $x_j^{\circ} \in N_{s_j}^{\circ}$ with $s_1 \neq s_2 \neq \cdots \neq s_n$. Similarly one has  
\begin{equation*}
E^{\psi}\left(\lambda_s\left(x\right)\right) = \lambda_s\left(E_s\left(x\right)\right) = 
\lambda\left(\iota_s^{-1}\circ E_s\left(x\right)\right), \quad x \in N_s. 
\end{equation*}
The conditional expectation $E^{\psi}$ can be shown to be independent from the choice of $\varphi$ (see the proposition below for more precise), and hence we rewrite $E := E^{\psi}$. The pair $\left(N, E\right)$ constructed so far is the desired one of von Neumann algebra and conditional expectation, and it is characterized by freeness with amalgamation as follows. 

\begin{fact} {\rm (\cite[\S\S5.6]{voiculescu:lecturenotes1132}; also see \cite[Proposition 2.5]{ueda:pacific}.)}\label{Fact2-2} 
Let $P$ be a von Neumann algebra with a normal $*$-isomorphism 
$\pi : D \rightarrow P$. Suppose that there are normal $*$-isomorphisms 
$\pi_s : N_s \rightarrow P$ with $\pi_s\circ\iota_s = \pi$ and a faithful normal conditional expectation $F : P \rightarrow \pi(D)$ such that 
\begin{itemize} 
\item the $\pi_s\left(N_s\right)$'s generate the whole $P${\rm ;}
\item $F\circ\pi_s = \pi\circ\iota_s^{-1}\circ E_s$ for every $s \in S${\rm ;}
\item the $\pi_s\left(N_s\right)$'s are free with amalgamation with respect to $F$. 
\end{itemize}
Then, there is a unique surjective normal $*$-isomorphism $\Pi : N \rightarrow P$ such 
that $\Pi\circ\lambda_s = \pi_s$ for every $s \in S$ and $\Pi\circ E = F\circ\Pi$.    
\end{fact} 

Since $\psi = \left(\psi\big|_{\lambda\left(D\right)}\right)\circ E = \varphi\circ\lambda^{-1}\circ E$, we see that $\sigma_t^{\psi}  = \underset{s \in S}{\bigstar_D} \sigma_t^{\varphi\circ\iota_s^{-1}\circ E_s}$ ($t \in {\mathbf R}$), where the right hand side is understood as free product of $*$-automorphisms (constructed based on the characterization by freeness, see e.g. \cite[p.366]{ueda:pacific}), i.e., 
\begin{equation*}
\left(\underset{s \in S}{\bigstar_D} \sigma_t^{\varphi\circ\iota_s^{-1}\circ E_s}\right)\left(\lambda_s(x)\right) :=\lambda_s\left(\sigma_t^{\varphi\circ\iota_s^{-1}\circ E_s}(x)\right), \quad x \in N_s.
\end{equation*}
Thanks to Connes' cocycle Radon-Nikodym theorem (see \cite[Vol.II; Chap.~VIII, \S3]{takesaki:book}), this formula of modular automorphisms is still valid even for every semifinite weight:  

\begin{proposition}\label{Prop2-3} {\rm (\cite[Theorem 2.6]{ueda:pacific})} For a faithful normal semifinite weight $\phi$ on $D$ we have 
\begin{equation*}
\sigma_t^{\phi\circ\lambda^{-1}\circ E} = \underset{s \in S}{\bigstar_D} \sigma_t^{\phi\circ\iota_s^{-1}\circ E_s} 
\quad (t \in {\mathbf R}). 
\end{equation*}
\end{proposition}

\section{Construction and Characterization} 

One would encounter ``difficulty" in dealing with conditional expectations (in connection with ``shift automorphisms") if straightforward adaptation of one of the group theoretic constructions of HNN extensions (see e.g.~\cite[Chap.~I, \S1.4]{serre:book}) was attempted in the von Neumann algebra setting. This forced us to seek for another route towards the construction of reduced HNN extensions. The rough idea is still essentially the same, but our method is completely different, avoiding the use of ``shift automorphisms" on ``infinite free products with amalgamations." The method is based on a simple fact on ``matrix multiplications" that we observed in our previous investigation on the reduced algebra of a certain amalgamated free product by a projection, see \cite[\S7]{ueda:us-jpn}.   

\medskip
Let $N$ be a $\sigma$-finite von Neumann algebra and $D$ be a distinguished von Neumann subalgebra with a faithful normal conditional expectation $E_D^N : N \rightarrow D$. Let us suppose that we have an (at most countably infinite) family $\Theta$ of normal $*$-isomorphisms $\theta : D \rightarrow  N$ with faithful normal conditional expectations $E_{\theta\left(D\right)}^N : N \rightarrow \theta\left(D\right)$. 

Set $\Theta_1 :=\{1 :=\mathrm{Id}_D\}\sqcup\Theta$, a disjoint union. 
Let us define the normal $*$-isomorphism $\iota_{\Theta} : D\otimes \ell^{\infty}\left(\Theta_1\right) \rightarrow N\otimes B\left(\ell^2\left(\Theta_1\right)\right)$ by 
\begin{equation*}
\iota_{\Theta}\left(x\otimes e_{\theta \theta}\right) := 
\begin{cases} 
x\otimes e_{1 1} & \text{if $\theta = 1$}, \\
\theta(x)\otimes e_{\theta \theta} & \text{if $\theta \in \Theta$}
\end{cases} 
\end{equation*}
for each $x \in D$, where the $e_{\theta_1 \theta_2}$'s denote the canonical matrix unit system 
in $B\left(\ell^2\left(\Theta_1\right)\right)$. Namely, in the operator matrix representation, we have 
\begin{equation*}
\iota_{\Theta} = 
\left[ 
\begin{array}{cccc} 
1 &  & & \\
 & \ddots & & \\
 & & \theta & \\
 & & &\ddots
\end{array}
\right], \quad 
\iota_{\Theta}\left(D\otimes\ell^{\infty}\left(\Theta_1\right)\right) = 
\left[ 
\begin{array}{cccc} 
D & & & \\
& \ddots & &  \\
& & \theta(D) & \\
& && \ddots
\end{array}
\right].
\end{equation*}
We also define the faithful normal conditional expectation 
$E_{\Theta} : N\otimes B\left(\ell^2\left(\Theta_1\right)\right) \rightarrow \iota_{\Theta}\left(D\otimes\ell^{\infty}\left(\Theta_1\right)\right)$
by 
\begin{equation*}
E_{\Theta} := 
\left[ 
\begin{array}{cccc} 
E_D^N &  &  & \\
 & \ddots & & \\
 & & E_{\theta\left(D\right)}^N & \\
 &  & &\ddots
\end{array}
\right]   
= 
\left(\sideset{}{^{\oplus}}\sum_{\theta \in \Theta_1} 
E_{\theta\left(D\right)}^N\otimes\mathrm {Id}_{{\mathbf C}e_{\theta \theta}}\right)\circ 
\left(\mathrm{Id}_N\otimes E_{\ell^{\infty}}\right), 
\end{equation*}
where $E_{\ell^{\infty}}$ is the unique conditional expectation from 
$B\left(\ell^2\left(\Theta_1\right)\right)$ onto $\ell^{\infty}\left(\Theta_1\right)$. 
Let us denote by $\iota_1$ the inclusion map of $D\otimes \ell^{\infty}\left(\Theta_1\right)$ into $N\otimes B\left(\ell^2\left(\Theta_1\right)\right)$, and define the faithful normal 
conditional expectation $E_1 : N\otimes B\left(\ell^2\left(\Theta_1\right)\right) \rightarrow D\otimes\ell^{\infty}\left(\Theta_1\right)$ by 
\begin{equation*}
E_1 := 
\left[ 
\begin{array}{cccc}
E_D^N &        &       &        \\
      & \ddots &       &        \\
      &        & E_D^N &        \\
      &        &       & \ddots
\end{array}
\right] = 
\left(E_D^N\otimes\mathrm{Id}_{\ell^{\infty}\left(\Theta_1\right)}\right)\circ\left(\mathrm{Id}_N
\otimes E_{\ell^{\infty}}\right) = \left(E_D^N\otimes E_{\ell^{\infty}}\right). 
\end{equation*}
We then construct the reduced free product with amalgamation: 
\begin{equation*}
\left({\mathcal N}, {\mathcal E}\right) = 
\left(N\otimes B\left(\ell^2\left(\Theta_1\right)\right), E_{\Theta} : \iota_{\Theta}\right) 
\underset{D\otimes\ell^{\infty}\left(\Theta_1\right)}{\bigstar} 
\left(N\otimes B\left(\ell^2\left(\Theta_1\right)\right), E_1 : \iota_1\right).     
\end{equation*}
The embedding maps of $N\otimes B\left(\ell^2\left(\Theta_1\right)\right)$ onto the 1st/2nd free components are denoted by $\lambda_{\Theta}$ and $\lambda_1$, respectively, and the embedding map of $D\otimes\ell^{\infty}\left(\Theta_1\right)$ into ${\mathcal N}$ by $\lambda$, i.e., $\lambda=\lambda_{\Theta}\circ\iota_{\Theta}=\lambda_1\circ\iota_1$. The desired HNN extension of $N$ by $\Theta$ with respect to $E_D^N$ and the $E_{\theta(D)}^N$'s will be constructed inside a corner subalgebra of ${\mathcal N}$. 

Let us define 
\begin{equation*}
u(\theta) := \lambda_1\left(e_{1 \theta}\right)\lambda_{\Theta}\left(e_{\theta 1}\right) 
\end{equation*}
with identifying $e_{\theta_1 \theta_2} = 1\otimes e_{\theta_1 \theta_2}$, and the following equation is a key to our construction:   
\begin{equation*} 
u(\theta)\lambda_{\Theta}\left(\theta(d)\otimes e_{1 1}\right) u(\theta)^* = \lambda_{\Theta}\left(d\otimes e_{1 1}\right), \quad d \in D, 
\end{equation*}
which simply comes from $\lambda_{\Theta}\left(\theta(d)\otimes e_{\theta\theta}\right) = \lambda\left(d\otimes e_{\theta\theta}\right) = \lambda_1\left(d\otimes e_{\theta\theta}\right)$ for each $d \in D, \theta \in \Theta_1$.
We also define the projection 
\begin{equation*}
p := \lambda\left(e_{11}\right) = \lambda_{\Theta}\left(e_{11}\right) 
\in {\mathcal N}, 
\end{equation*}
and then introduce the unital normal $*$-isomorphism $\pi$ from $N$ into the corner subalgebra $p{\mathcal N}p$ defined by
\begin{equation*}
\pi(n) := \lambda_{\Theta}\left(n\otimes e_{11}\right), \quad n \in N. 
\end{equation*}
The partial isometries $u(\theta)$'s, can be thought of as unitaries in the corner $p{\mathcal N}p$ since their left and right supports are the projection $p$, and the above-mentioned key equation is translated into the following algebraic relation: 
\begin{equation}\label{Eq2}
u(\theta)\pi\left(\theta\left(d\right)\right)u(\theta)^* = \pi\left(d\right), \quad d \in D. 
\end{equation}
Set 
\begin{equation*}
M := \pi\left(N\right) \vee \left\{ u(\theta) : \theta \in \Theta \right\}'' \subseteq p{\mathcal N}p. 
\end{equation*}
Let us consider a faithful normal semifinite weight on $D\otimes\ell^{\infty}\left(\Theta_1\right)$: 
\begin{equation*}
\phi :=\varphi\otimes
\left(\mathrm{Tr}\big|_{\ell^{\infty}\left(\Theta_1\right)}\right)
\end{equation*}
with a faithful normal state $\varphi$ on $D$, where $\mathrm{Tr}$ is the (non-normalized) canonical normal trace on $B\left(\ell^2\left(\Theta_1\right)\right)$.   

\begin{lemma}\label{Lem3-1} We have 
\begin{eqnarray}
\sigma_t^{\phi\circ\lambda^{-1}\circ{\mathcal E}}\left(p\right) &=& p, \label{Eq3}\\ 
\sigma_t^{\phi\circ\lambda^{-1}\circ{\mathcal E}}\left(\pi\left(n\right)\right) &=& 
\pi\left(\sigma_t^{\varphi\circ E_D^N}\left(n\right)\right), \label{Eq4}\\ 
\sigma_t^{\phi\circ\lambda^{-1}\circ{\mathcal E}}\left(u(\theta)\right) &=& u(\theta)
\pi\left(
\left[D\varphi\circ\theta^{-1}\circ E_{\theta(D)}^N : D\varphi\circ E_D^N\right]_t
\right) \label{Eq5}
\end{eqnarray}
for each $t \in {\mathbf R}$, $n \in N$, $\theta \in \Theta$. 
\end{lemma} 
\begin{proof} The equations \eqref{Eq3}, \eqref{Eq4} are straightforward from Proposition \ref{Prop2-3}, while the last one needs some additional efforts. In fact, we have 
\begin{align*}
\sigma_t^{\phi\circ\lambda^{-1}\circ{\mathcal E}}\left(u(\theta)\right) 
&= 
\lambda_1
\left(
\sigma_t^{\phi\circ\iota_1^{-1}\circ E_1}\left(1\otimes e_{1 \theta}\right)\right) 
\lambda_{\Theta}\left(\sigma_t^{\phi\circ\iota_{\Theta}^{-1}\circ E_{\Theta}}
\left(1\otimes e_{\theta 1}\right)
\right) \\ 
&= 
\lambda_1\left(1\otimes e_{1 \theta}\right) 
\lambda_{\Theta}\left(
\left[D\varphi\circ\theta^{-1}\circ E_{\theta(D)}^N : D\varphi\circ E_D^N\right]_t\otimes 
e_{\theta 1}\right) \\
&= 
\lambda_1\left(e_{1 \theta}\right) 
\lambda_{\Theta}\left(e_{\theta 1}\right)
\lambda_{\Theta}\left(
\left[D\varphi\circ\theta^{-1}\circ E_{\theta(D)}^N : D\varphi\circ E_D^N\right]_t\otimes 
e_{1 1}\right) \\
&= 
u(\theta) \pi\left(
\left[D\varphi\circ\theta^{-1}\circ E_{\theta(D)}^N : D\varphi\circ E_D^N\right]_t
\right), 
\end{align*}
where the second equality comes from the so-called ``balanced weight technique" due to Connes (see \cite[Vol.II; Chap.~VIII, \S3, p.111--113]{takesaki:book}).  
\end{proof} 

Since $\pi(N) = p\lambda_{\Theta}\left(N\otimes B\left(\ell^2\left(\Theta_1\right)\right)\right)p$, the restriction of the normal conditional expectation 
\begin{equation*}
{\mathcal E}_{\Theta} : {\mathcal N} \rightarrow 
\lambda_{\Theta}\left(N\otimes B\left(\ell^2\left(\Theta_1\right)\right)\right)
\end{equation*} 
that preserves $\phi\circ\lambda^{-1}\circ{\mathcal E}$ (and hence ${\mathcal E} = {\mathcal E}\circ{\mathcal E}_{\Theta}$ holds) to $M$ gives a faithful normal conditional expectation  
\begin{equation*}
E^M_{\pi(N)} :={\mathcal E}_{\Theta}\big|_M : M \rightarrow \pi(N).  
\end{equation*}
We have $\left(\phi\circ\lambda^{-1}\circ{\mathcal E}\right)\big|_M = \varphi\circ E_D^N\circ\pi^{-1}\circ E_{\pi(N)}^M$, and hence, by Takesaki's theorem \cite[Vol.II; Theorem 1.2 in Chap.~VIII, \S1]{takesaki:book} we get  
\begin{equation}\label{Eq6}
\sigma_t^{\varphi\circ E_D^N\circ\pi^{-1}\circ E_{\pi(N)}^M} = \sigma_t^{\phi\circ\lambda^{-1}\circ{\mathcal E}}\big|_M, \quad t \in {\mathbf R} 
\end{equation}  
since $\sigma_t^{\phi\circ\lambda^{-1}\circ{\mathcal E}}(M) = M$ for every $t \in {\mathbf R}$ thanks to Lemma \ref{Lem3-1}. 

\begin{definition}\label{Def3-1} {\rm {\bf (Reduced HNN extensions)}} We call the pair $\left(M, E^M_{\pi(N)}\right)$ constructed so far the reduced HNN extension {\rm (}or HNN extension, in short{\rm )} of $N$ by $\Theta$ with respect to $E_D^N$ and the $E_{\theta(D)}^N$'s, and denote it by    
\begin{equation*}
\left(M, E^M_{\pi(N)}\right) = \left(N, E_D^N\right) \underset{D}{\bigstar} 
\left(\Theta, \left\{E_{\theta(D)}^N\right\}_{\theta \in \Theta}
\right). 
\end{equation*} 
When no confusion occurs, we will write $M = N \underset{D}{\bigstar} \Theta$ for short. The given von Neumann algebra $N$ is called the base algebra, and each $u(\theta)$ the stable unitary of $\theta \in \Theta$. 
\end{definition}  

\begin{definition}\label{Def3-2} {\rm {\bf (Reduced words)}} An element {\rm (}in $M${\rm )} 
\begin{equation*}
w =
u\left(\theta_0\right)^{\varepsilon_0}\pi\left(n_1\right)u\left(\theta_1\right)^{\varepsilon_1}\pi\left(n_2\right)
\cdots\pi\left(n_{\ell}\right)u\left(\theta_{\ell}\right)^{\varepsilon_{\ell}}
\end{equation*} 
with $n_1, n_2,\dots, n_{\ell} \in N$, $\theta_0, \theta_1,\dots, \theta_{\ell} \in \Theta$, $\varepsilon_0, \varepsilon_1, \dots, \varepsilon_{\ell} \in \left\{1, -1\right\}$ {\rm (}possibly with $w = u\left(\theta_0\right)^{\varepsilon_0}${\rm )} is called a reduced word {\rm (}or said to be of reduced form{\rm )} if $\theta_{j-1} = \theta_{j}$ and 
$\varepsilon_{j-1} \neq \varepsilon_j$ imply that  
\begin{itemize}
\item $n_j \in N_{\theta}^{\circ} :=
\mathrm{Ker}E_{\theta(D)}^N$ with $\theta :=
\theta_{j-1} = \theta_j$, when $\varepsilon_{j-1}  = 1, \ \varepsilon_j = -1${\rm ;} 
\item $n_j \in N^{\circ} :=\mathrm{Ker}E_D^N$,  
when $\varepsilon_{j-1}  = -1, \ \varepsilon_j = 1$.  
\end{itemize}
\end{definition} 

We should point out that our definition of reduced words agrees with so-called Britton's lemma in combinatorial group theory (see \cite[p.181]{lyndonschupp:book}), where a reduced word is named as a normal form, and the sets of representatives of right cosets of distinguished subgroups should be regarded as the counterparts of $N^{\circ}$ and the $N_{\theta}^{\circ}$'s in our consideration. 

\begin{remark}\label{Rem3-2} {\rm It is plain to see that 
$u\left(\theta_1\right)^{\varepsilon_1} \cdots u\left(\theta_{\ell}\right)^{\varepsilon_{\ell}}$ is of reduced form 
in the above sense if and only if so is  $\theta_1^{\varepsilon_1}\cdots\theta_{\ell}^{\varepsilon_{\ell}}$ 
in the free group ${\mathbb F}\left(\Theta\right)$ over the generating set $\Theta$.}  
\end{remark}

\begin{definition}\label{Def3-3} {\bf (Conditions needed for characterization)} We introduce the following two conditions: 
\begin{itemize}
\item[(A)] $u(\theta)\pi\left(\theta\left(d\right)\right)u(\theta)^* = \pi\left(d\right)$ 
for every $d \in D$, $\theta \in \Theta$. 
\item[(M)] For every reduced word $w$, one has $E^M_{\pi(N)}\left(w\right) = 0$. 
\end{itemize}
\end{definition} 

\begin{theorem}\label{Thm3-3} 
The pair $\left(M, E^M_{\pi(N)}\right)$ constructed above satisfies the conditions {\rm (A)}, {\rm (M)}.  On the other hand, the conditions {\rm (A)}, {\rm (M)} characterize the pair $\left(M, E^M_{\pi(N)}\right)$ completely under the assumption that $\pi(N)$ and the $u(\theta)$'s generate $M$ as von Neumann algebra. Strictly speaking, the conditional expectation of the pair in question is completely determined by those conditions. 
\end{theorem} 
\begin{proof} Let us denote the 1st/2nd free components of ${\mathcal N}$ by ${\mathcal N}_{\Theta}$, ${\mathcal N}_1$, respectively, for short, i.e., ${\mathcal N}_{\Theta} := \lambda_{\Theta}\left(N\otimes B\left(\ell^2\left(\Theta_1\right)\right)\right)$, ${\mathcal N}_1 := \lambda_1\left(N\otimes B\left(\ell^2\left(\Theta_1\right)\right)\right)$, 
and set ${\mathcal N}_{\Theta}^{\circ} :={\mathcal N}_{\Theta} \cap\mathrm{Ker}{\mathcal E} = \lambda_{\Theta}\left(\mathrm{Ker}E_{\Theta}\right)$ and ${\mathcal N}_1^{\circ} :={\mathcal N}_1\cap\mathrm{Ker}{\mathcal E} = \lambda_1\left(\mathrm{Ker}E_1\right)$ as usual. 

The condition (A) was already verified, see the equation \eqref{Eq2}, and thus it suffices to check the condition (M) for the first half of the assertions. Let us choose a word
\begin{equation*}
w = u\left(\theta_0\right)^{\varepsilon_0}\pi\left(n_1\right)u\left(\theta_1\right)^{\varepsilon_1}\pi\left(n_2\right)
\cdots\pi\left(n_{\ell}\right)u\left(\theta_{\ell}\right)^{\varepsilon_{\ell}}, 
\end{equation*}
and then we have 
{\allowdisplaybreaks 
\begin{align*}
w 
&= 
\left(
\lambda_1\left(e_{1 \theta_0}\right)\lambda_{\Theta}\left(e_{\theta_0 1}\right)
\right)^{\varepsilon_0} 
\lambda_{\Theta}\left(n_1\otimes e_{11}\right) 
\left(
\lambda_1\left(e_{1 \theta_1}\right)\lambda_{\Theta}\left(e_{\theta_1 1}\right)
\right)^{\varepsilon_1} \\
&\phantom{aaaaaaaaaaaaaaaaaaaaaaaaaaaaaa}\cdots 
\lambda_{\Theta}\left(n_{\ell}\otimes e_{11}\right)
\left(
\lambda_1\left(e_{1 \theta_{\ell}}\right)\lambda_{\Theta}\left(e_{\theta_{\ell} 1}\right)
\right)^{\varepsilon_0}.   
\end{align*}
}Here, we briefly explain how to manipulate this word in a typical case: If $\varepsilon_{j-1} = -1$, $\varepsilon_j = 1$, then 
{\allowdisplaybreaks 
\begin{align*}
&\left(\lambda_1\left(e_{1 \theta_{j-1}}\right)\lambda_{\Theta}\left(e_{\theta_{j-1} 1}\right)
\right)^{\varepsilon_{j-1}} 
\lambda_{\Theta}\left(n_j\otimes e_{11}\right) 
\left(\lambda_1\left(e_{1 \theta_j}\right)\lambda_{\Theta}\left(e_{\theta_j 1}\right)
\right)^{\varepsilon_j} \\
&= 
\lambda_{\Theta}\left(e_{1 \theta_{j-1}}\right)
\lambda_1\left(E_D^N\left(n_j\right)\otimes e_{\theta_{j-1} \theta_j}\right)
\lambda_{\Theta}\left(e_{\theta_j 1}\right) \\
&\phantom{aaaaaaa}+
\lambda_{\Theta}\left(e_{1 \theta_{j-1}}\right)
\lambda_1\left(e_{\theta_{j-1} 1}\right)
\lambda_{\Theta}\left(n_j^{\circ}\otimes e_{11}\right)
\lambda_1\left(e_{1 \theta_j}\right)
\lambda_{\Theta}\left(e_{\theta_j 1}\right)
\end{align*}
}with $n_j^{\circ} = n_j - E_D^N\left(n_j\right)$. If $n_j \in N^{\circ}$, then this belongs to ${\mathcal N}_{\Theta}^{\circ}{\mathcal N}_1^{\circ}{\mathcal N}_{\Theta}^{\circ}{\mathcal N}_1^{\circ}{\mathcal N}_{\Theta}^{\circ}$ since the first term disappears in this case. On the other hand, if $n_j$ is arbitrary but $\theta_{j-1} \neq \theta_j$, then it belongs to ${\mathcal N}_{\Theta}^{\circ}{\mathcal N}_1^{\circ}{\mathcal N}_{\Theta}^{\circ} + {\mathcal N}_{\Theta}^{\circ}{\mathcal N}_1^{\circ}{\mathcal N}_{\Theta}^{\circ}{\mathcal N}_1^{\circ}{\mathcal N}_{\Theta}^{\circ}$. 
In this way, one can easily observes that, if the word $w$ is of reduced form, then it belongs to the linear span of alternating words in ${\mathcal N}_{\Theta}^{\circ}$ and ${\mathcal N}_1^{\circ}$ of length  greater than $2$. Therefore, we have ${\mathcal E}_{\Theta}\left(w\right) = 0$, which asserts the condition (M). 

Next, we will show the latter half of the assertions. To do so, it is enough to explain how one can compute the moment-value: 
\begin{equation*}
E^M_{\pi(N)}\left(u\left(\theta_0\right)^{\delta_0}\pi\left(x_1\right)u\left(\theta_1\right)^{\delta_1}\pi\left(x_2\right)\cdots\pi\left(x_m\right)u\left(\theta_m\right)^{\delta_m}\right) 
\end{equation*}
of any given $x_1, x_2,\dots, x_m \in N$, $\theta_0, \theta_1, \dots, \theta_m \in \Theta$, $\delta_0, \delta_1, \dots, \delta_m \in {\mathbb Z}\setminus\{0\}$, by using only the conditions (A), (M). In fact, if the resulting value could be expressed uniquely in terms of only the data of $\left(\theta_0, \delta_0\right), x_1, \left(\theta_1, \delta_1\right),\dots, x_m, \left(\theta_m, \delta_m\right)$ together with $E^N_D$ and the $E^N_{\theta(D)}$'s, then the desired assertion would follow. Our technique is the essentially same as in the case of free products with amalgamations. Namely, we use the decompositions: 
\begin{equation*}
n = E_D^N\left(n\right) + n^{\circ}\ \text{or}\ 
E_{\theta(D)}^N\left(n\right) + [n]_{\theta}^{\circ}\ \ n \in N,
\end{equation*}
where we define $[n]_{\theta}^{\circ} :=n-E_{\theta(D)}^N\left(n\right)$. By the repeated use of the decompositions together with the condition (A), we can make the moment-value in question a (finite) sum of the form: 
\begin{equation*}
\sum_{\text{$w$: reduced word or $1$}} \pi\left(n(w)\right) E_{\pi(N)}^M(w)
\end{equation*}
with coefficients $n(w)$ being words in $D$ and $\theta(D)$ (in $N$), and all the coefficients $n(w)$ and all the words $w$ (the moment-value $E^M_{\pi(N)}(w)$ takes $0$ if $w$ is of reduced form or otherwise, $w = 1$) appeared in the above expression are uniquely determined from the given data $\left(\theta_0,\delta_0\right)$, $x_1$, $\left(\theta_1, \delta_1\right)$, $\dots$, $x_m$, $\left(\theta_m, \delta_m\right)$ together with $E^N_D$ and the $E^N_{\theta(D)}$'s. Therefore, our desired assertion follows. 
\end{proof} 

Let $u\left(g\right)$, $g \in {\mathbb F}\left(\Theta\right)$, be the natural group isomorphism from the free group ${\mathbb F}\left(\Theta\right)$ into the unitary group ${\mathcal U}\left(M\right)$ given by the correspondence $\theta \in \Theta \longmapsto u(\theta) \in {\mathcal U}\left(M\right)$. Let us denote by $\ell(\ \cdot\ )$ the usual word length function with respect to the generating set $\Theta$. The computation given in the above proof implies the following corollary: 

\begin{corollary}\label{Cor3-4} Let 
$w = u\left(\theta_0\right)^{\varepsilon_0}\pi\left(n_1\right)u\left(\theta_1\right)^{\varepsilon_1}\pi\left(n_2\right)\cdots\pi\left(n_{\ell}\right)u\left(\theta_{\ell}\right)^{\varepsilon_{\ell}}$ be a word in $M$, and set 
$g := \theta_0^{\varepsilon_0}\theta_1^{\varepsilon_1}\cdots\theta_{\ell}^{\varepsilon_{\ell}}$, 
a word in ${\mathbb F}\left(\Theta\right)$ {\rm (}obtained by replacing all $n_j$'s by the identity $1${\rm )}. Then we have   
\begin{equation*}
\ell\left(g\right) \neq 0 \Longrightarrow E_{\pi(N)}^M\left(w\right) = 0.
\end{equation*}
In particular, the unitaries $u(\theta)$'s form a free family of Haar unitaries, so that they generate the free group factor $L\left({\mathbb F}\left(\Theta\right)\right)$. 
\end{corollary} 

The following corollary is also straightforward from Theorem \ref{Thm3-3}: 

\begin{corollary}\label{Cor3-5} Let $G *_H \theta = \langle G, t : t\theta(h)t^{-1} = h,\ h \in H \rangle$ be an HNN extension of base group $G$ with stable letter $t$ by group isomorphism $\theta$ from $H$ into $G$. Then, the group von Neumann algebra $L\left(G *_H \theta\right)$ can be identified with the reduced HNN extension of the  base algebra $L\left(G\right)$ with the stable unitary $\lambda\left(t\right)$, where all the necessary conditional expectations are chosen as the canonical tracial state preserving ones. 
\end{corollary}

We then discuss what phenomenon occurs when $D$ and the $\theta(D)$'s are assumed to be all mutually inner conjugate. Let $\left(M,E_{N}^M\right)$ be as above with identifying $n = \pi(n)$, $n \in N$. We here suppose that every $\theta \in \Theta$ has a unitary $w_{\theta} \in N$ with the following properties: (i) $\mathrm{Ad}w_{\theta}\circ\theta \in \mathrm{Aut}(D)$; and (ii) $E_{\theta(D)}^N = \mathrm{Ad}w_{\theta}^*\circ E_D^N\circ\mathrm{Ad}w_{\theta}$. Define the action $\gamma$ of $\mathbb{F}(\Theta)$ on $D$ in such a way that $\gamma_{\theta} = \mathrm{Ad}w_{\theta}\circ\theta$, and consider the free product with amalgamation: 
\begin{equation*}
\left(L, F_D^L\right) := 
\left(N,E_D^N\right) \underset{D}{\bigstar}
\left(D\rtimes_{\gamma}\mathbb{F}(\Theta), E_D^{D\rtimes_{\gamma}\mathbb{F}(\Theta)}\right), 
\end{equation*}
where $E_D^{D\rtimes_{\gamma}\mathbb{F}(\Theta)}$ is the canonical conditional expectation. The faithful normal conditional expectation from $L$ onto the 1st free component $N$ that preserves $F_D^L$ is denoted by $F_N^L$. Then, we have the following simple corollary: 
 
\begin{corollary}\label{Cor3-6} 
In the above setting, the correspondence{\rm :} 
\begin{align*}
n \in M &\mapsto n \in L,\quad n \in N; \\
u(\theta) \in M &\mapsto \lambda^{\gamma}(\theta)^* w_{\theta} \in L,\quad \theta \in \Theta
\end{align*}
gives a $*$-isomorphism between $M$ and $L$ that intertwines $E_N^M$ and $F_N^L$. Here, $\lambda^{\gamma} : \mathbb{F}(\Theta) \rightarrow D\rtimes_{\gamma}\mathbb{F}(\Theta)$ {\rm (}$\subseteq L${\rm )} denotes the canonical unitary representation. 
\end{corollary} 
\begin{proof} It is plain to verify that the pair $\left(L \supseteq N, F^L_N\right)$ with the unitaries $\lambda^{\gamma}\left(\theta\right)^*w_{\theta}$, $\theta \in \Theta$, satisfies the conditions (A), (M) with respect to $\Theta$ and $E_D^N$, $\left\{E_{\theta(D)}^N\right\}_{\theta\in\Theta}$. In fact, the condition (A) follows from the above (i), while the (M) from the fact that $x \in N_{\theta}^{\circ} = \mathrm{Ker}E^N_{\theta(D)}$ if and only if $w_{\theta}x w_{\theta}^* \in N^{\circ} = \mathrm{Ker}E_D^N$ thanks to the above (ii).  
\end{proof}     

\begin{remarks}\label{Rem3-7} {\rm (1) {\bf [HNN extensions arising from inner conjugate Cartan subalgebras]} Assume that $N$ is a non-type I factor with separable predual (or more generally, a von Neumann algebra with separable predual having no type I direct summand) and further that $D$ and the $\theta(D)$'s are all Cartan subalgebras in $N$. By the uniqueness of normal conditional expectations onto those Cartan subalgebras, if those Cartan subalgebras are all mutually inner conjugate, then Corollary \ref{Cor3-6} enables us to apply our previous results \cite{ueda:pacific}\cite{ueda:us-jpn}\cite{ueda:transAMS} to the HNN extension $M = N \bigstar_D \Theta$ without any change. However, we have no general result without this inner conjugacy assumption among Cartan subalgebras in question. 

(2) A special case of Corollary \ref{Cor3-6} was one of the starting points of the present work. In fact, in the setting of Corollary \ref{Cor3-6}, the group theoretic construction based on shift automorphisms on infinite amalgamated free products is valid when all the $w_{\theta} = 1$ (so that $\gamma_{\theta} = \theta$). Concerning this, we point out that the amalgamated free product appeared in Corollary \ref{Cor3-6} has the crossed-product decomposition{\rm :} 
\begin{equation*}
M = N\left(\Theta\right)\rtimes{\mathbb F}\left(\Theta\right)
\end{equation*}
by the free Bernoulli shift on  
\begin{equation*}
N\left(\Theta\right) :=
\underset{g \in {\mathbb F}\left(\Theta\right)}{\bigstar_D}
\left(N, E_D^N : \gamma_g\right). 
\end{equation*}
 (See e.g. \cite[\S3]{hiaiueda:jot}, where only the case of $D = {\mathbf C}1$ was treated, but the argument works even in this case.) More on this will be discussed in the next section with full generality. }
\end{remarks}


\section{Modular Theory} 

Let 
\begin{equation*}
\left(M, E^M_{\pi(N)}\right) = \left(N, E_D^N\right) \underset{D}{\bigstar} 
\left(\Theta, \left\{E_{\theta(D)}^N\right\}_{\theta \in \Theta}
\right)  
\end{equation*}
be the HNN extension of base von Neumann algebra $N$ with stable unitaries $u(\theta)$, $\theta \in \Theta$. Here, we will use the construction and the notations of HNN extensions given in the previous section; however, in what follows, we will identify $n = \pi\left(n\right)$, $n \in N$, so $\pi$ will be omitted. The next theorem is immediately derived from Lemma \ref{Lem3-1} with the aid of Connes' cocycle Radon-Nikodym theorem (see \cite[Vol.II; Chap.~VIII, \S3]{takesaki:book}).

\begin{theorem} \label{Thm4-1} For a faithful normal semifinite weight $\psi$ on $D$, we have 
\begin{equation*}
\sigma_t^{\psi\circ E_D^N\circ E^M_N}\left(u(\theta)\right) = 
u(\theta) 
\left[D\psi\circ\theta^{-1}\circ E_{\theta(D)}^N : D\psi\circ E_D^N\right]_t, 
\quad 
t \in {\mathbf R}. 
\end{equation*}
\end{theorem}

This theorem implies the following criterion for the existence of traces on HNN extensions:  

\begin{corollary}\label{Cor4-2} If $N$ has a faithful normal semifinite trace $\tau$ and if the given $E_D^N$ and the $E_{\theta(D)}^N$'s satisfy the relation{\rm :}
\begin{equation*}
\tau = 
\left(\tau\big|_D\right)\circ E_D^N = 
\left(\tau\big|_D\right)\circ\theta^{-1}\circ E_{\theta(D)}^N, \quad 
\theta \in \Theta, 
\end{equation*}
then so does $M$, and more precisely $\tau\circ  E^M_N$ becomes a trace. In particular, if $N$ is semifinite with a faithful normal semifinite trace $\tau$ and if the given conditional expectations are $\tau$-preserving and $\tau|_{\theta(D)} = \left(\tau|_D\right)\circ\theta$ holds, then $\tau\circ E_N^M$ is a trace. 
\end{corollary}   

A crossed-product decomposition fact for HNN extensions was given in Remarks \ref{Rem3-7}, (2) under a very special assumption. Here, we give such a fact with full generality. 

\begin{corollary}\label{Cor4-3} Let us denote by $N\left(\Theta\right)$ the von Neumann subalgebra generated by all the $u\left(g\right) N u\left(g\right)^*$, $g \in {\mathbb F}\left(\Theta\right)$. Then we have the crossed-product decomposition of $M${\rm :}
\begin{equation*}
M = N\left(\Theta\right)\rtimes_{\mathrm{Ad}u} {\mathbb F}\left(\Theta\right) 
\end{equation*}
with the natural adjoint action $\mathrm{Ad}u: g \in {\mathbb F}\left(\Theta\right) \mapsto 
\mathrm{Ad}u\left(g\right) \in \mathrm{Aut}\left(N\left(\Theta\right)\right)$.  
\end{corollary}
\begin{proof} First of all, we should remark that Theorem \ref{Thm4-1} shows that there is a unique faithful normal conditional expectation from $M$ onto $N\left(\Theta\right)$ that preserves $E^M_N$. Thus, the desired assertion is derived from Corollary \ref{Cor3-4} together with the well-known characterization of discrete crossed-products in terms of conditional expectations. 
\end{proof} 

\begin{remark}\label{Rem4-4} {\rm Theorem \ref{Thm4-1} says that each subalgebra $u\left(g\right) N u\left(g\right)^*$, $g \in {\mathbb F}\left(\Theta\right)$ with $g \neq e$ ($e$ denotes the identity), is not necessary to be globally invariant under the modular automorphism associated with $\psi\circ E^M_N$.} 
\end{remark} 

Theorem \ref{Thm4-1} enables us to show that the continuous core of the HNN extension $M$ in question becomes again an HNN extension. For a better description, it is convenient to use a recent formulation of continuous cores due to S.~Yamagami \cite{yamagami:aspect-modular}. (See also A.~J.~Falcone and M.~Takesaki \cite{falconetakesaki:quantumflow}, and the reader may consult \cite[Vol.II; Chap.~XII, \S6]{takesaki:book} for more detailed account.) Following the formulation, the continuous core $\widetilde{P}$ of a given $P$ can be understood as an abstract von Neumann 
algebra generated by two kinds of symbols $x \in P$ and $\psi^{it}$ with a faithful normal semi-finite weight $\psi$ on $P$, which satisfy  the relations: 
\begin{equation*}
\psi^{it}x\psi^{-it} = \sigma_t^{\psi}\left(x\right), \quad \psi^{it}\psi^{is} = \psi^{i(t+s)}, \quad 
\phi^{it}\psi^{-it} = \left[D\phi: D\psi\right]_t
\end{equation*}
for faithful normal semi-finite weights $\phi$, $\psi$ on $P$. 
It is known that such a von Neumann algebra $\widetilde{P}$ can be realized as the crossed-product $P\rtimes_{\sigma^{\psi}}{\mathbf R}$, where $\psi^{it}$ denotes the canonical unitary implementation $\lambda^{\psi}\left(t\right)$ of ${\mathbf R}$ inside the crossed-product, and $\phi^{it} = \left[D\phi: D\psi\right]_t \psi^{it} = \left[D\phi: D\psi\right]_t \lambda^{\psi}\left(t\right)$ in general. 
 
In our setting, the inclusion relations $M \supseteq N \supseteq D$, $M \supseteq N \supseteq \theta(D)$, $\theta \in \Theta$, with the faithful normal conditional expectations $E^M_N : M \rightarrow N$, $E_D^N : \rightarrow D$, $E_{\theta(D)}^N : N \rightarrow \theta(D)$ give us the following natural embeddings and mapping: 
\begin{eqnarray*} 
\widetilde{D} \hookrightarrow \widetilde{N} \quad &\text{by}& \quad 
\begin{cases}
d \in D \overset{\mathrm{identify}}{\leftrightarrow} 
d \in D \subseteq N \subseteq \widetilde{N}, \\
\varphi^{it} \in \widetilde{D} \overset{\mathrm{identify}}{\leftrightarrow} 
\left(\varphi\circ E_D^N\right)^{it} \in \widetilde{N}; 
\end{cases} \\
\widetilde{\theta}: 
\widetilde{D} \rightarrow \widetilde{N} \quad &\text{by}& \quad 
\begin{cases}
d \in D \mapsto \theta\left(d\right) \in \theta(D) \subseteq N \subseteq \widetilde{N}, \\ 
\varphi^{it} \in \widetilde{D} \mapsto \left(\varphi\circ\theta^{-1}\circ E_{\theta(D)}^N\right)^{it} 
\in \widetilde{N}; 
\end{cases} \\
\widetilde{N} \hookrightarrow \widetilde{M} \quad &\text{by}& \quad 
\begin{cases} 
n \in N \overset{\mathrm{identify}}{\leftrightarrow} n \in M \subseteq \widetilde{M}, \\
\phi^{it} \in \widetilde{N} \overset{\mathrm{identify}}{\leftrightarrow}  \left(\phi\circ E^M_N\right)^{it} \in \widetilde{M}, 
\end{cases} 
\end{eqnarray*}
and the conditional expectations 
\begin{equation*}
\widehat{E_D^N} : \widetilde{N} \rightarrow \widetilde{D}, 
\quad
\widehat{E_{\theta(D)}^N} : \widetilde{N} \rightarrow 
\widetilde{\theta(D)} = \widetilde{\theta}\left(\widetilde{D}\right), 
\quad
\widehat{E^M_N} : \widetilde{M} \rightarrow \widetilde{N}
\end{equation*} 
constructed in such a way that 
\begin{eqnarray*}
\widehat{E_D^N}\big|_N = E_D^N, &\quad&
\widehat{E_D^N} \left(\left(\varphi\circ E_D^N\right)^{it}\right) = 
\left(\varphi\circ E_D^N\right)^{it}; \\ 
\widehat{E_{\theta(D)}^N}\big|_N = E_{\theta(D)}^N, &\quad& 
\widehat{E_{\theta(D)}^N} \left(\left(\varphi\circ\theta^{-1}\circ E_{\theta(D)}^N\right)^{it}\right) = 
\left(\varphi\circ\theta^{-1}\circ E_{\theta(D)}^N\right)^{it}; \\
\widehat{E^M_N}\big|_M = E^M_N, &\quad&
\widehat{E^M_N} \left(\left(\phi\circ E^M_N\right)^{it}\right) = 
\left(\phi\circ E^M_N\right)^{it}
\end{eqnarray*}
for faithful normal positive linear functionals $\varphi \in D_*$, $\phi \in N_*$, where one should remind the following formula:  
\begin{equation*}
\widetilde{\theta}\left(\left(\varphi\circ E_D^N\right)^{it}\right) \left(\ = \widetilde{\theta}\left(\varphi^{it}\right)\ \right) = 
\left(\varphi\circ\theta^{-1}\circ E_{\theta(D)}^N\right)^{it}.
\end{equation*}
For a faithful normal state $\varphi$ on $D$, we have, 
in the continuous core $\widetilde{M}$, 
\begin{equation}\label{Eq7}
\left(\varphi\circ E_D^N\circ E^M_N\right)^{it} u(\theta) =  u(\theta)
\left(\varphi\circ\theta^{-1}\circ E_{\theta(D)}^N\circ E^M_N\right)^{it}, \quad t \in \mathbf{R}, 
\end{equation} 
thanks to Theorem \ref{Thm4-1}. The general assertion given below is a simple application of the formula \eqref{Eq7} and Theorem \ref{Thm3-3}, i.e., the characterization of HNN extensions. 

\begin{theorem}\label{Thm4-5} The pair $\left(\widetilde{M},\widehat{E^M_N}\right)$ is again the HNN extension of the base algebra $\widetilde{N}$ by the family $\widetilde{\Theta}  :=\left\{ \widetilde{\theta}\ :\ \theta \in \Theta\right\}$ 
with the stable unitaries $u(\theta)$, $\theta \in \Theta$, with respect to the conditional expectations $\widehat{E_D^N}$ and $\widehat{E_{\theta(D)}^N}$, $\theta \in \Theta$, that is, 
\begin{equation*}
\left(\widetilde{M}, \widehat{E^M_N}\right) = 
\left(\widetilde{N}, \widehat{E_D^N}\right) \underset{\widetilde{D}}{\bigstar} 
\left(\widetilde{\Theta}, \left\{\widehat{E_{\theta(D)}^N}\right\}_{\theta \in \Theta}
\right).  
\end{equation*}
We will often denote this identification by $\widetilde{M} = \widetilde{N}\underset{\widetilde{D}}{\bigstar}\widetilde{\Theta}$ for short.  
\end{theorem}  
\begin{proof} The condition (A) follows from the formula \eqref{Eq7}. Indeed, for each $d \in D$, $t \in {\mathbf R}$, the formula \eqref{Eq7} enables us to compute   
\begin{equation*}
u(\theta)\widetilde{\theta}\left(d\left(\varphi\circ E_D^N\right)^{it}\right)u(\theta)^* = d \left(\varphi\circ E_D^N\circ E^M_N\right)^{it} = d \left(\varphi\circ E_D^N\right)^{it}. 
\end{equation*}
Hence, it suffices to verify the condition (M). The argument to do so is similar to that in \cite[Theorem 5.1]{ueda:pacific}. Let $\widetilde{w}$ be a word in $\widetilde{M}$, i.e., 
\begin{equation*}
\widetilde{w} = 
u\left(\theta_0\right)^{\varepsilon_0} \widetilde{n}_1 u\left(\theta_1\right)^{\varepsilon_1} \cdots \widetilde{n}_{\ell} u\left(\theta_{\ell}\right)^{\varepsilon_{\ell}}. 
\end{equation*}
with $\widetilde{n}_1,\dots, \widetilde{n}_{\ell} \in \widetilde{N}$. Then, Kaplansky's density theorem enables us to reduce our consideration to the case that each $\widetilde{n}_j$ is in a (not necessary common) suitable dense $*$-subalgebra of $\widetilde{N}$. Such a dense $*$-subalgebra is chosen as the $*$-algebra generated by $N$ and the $\left(\varphi\circ E_D^N\right)^{it}$, $t \in {\mathbf R}$, or by $N$ and the $\left(\varphi\circ\theta^{-1}\circ E_{\theta(D)}^N\right)^{it}$, $t \in {\mathbf R}$, and we can assume that each $\widetilde{n}_j$ is of the form: $n\left(\varphi\circ E_D^N\right)^{it}$ or $n\left(\varphi\circ\theta^{-1}\circ E_{\theta(D)}^N\right)^{it}$ thanks to \cite[Lemma 5.2]{ueda:pacific}. If $\widetilde{w}$ is of reduced form, then the repeated use of the formula \eqref{Eq7} enables us to make $\widetilde{w}$ a word of the form: 
{\allowdisplaybreaks 
\begin{align*}
\left(\text{a reduced word in $M$}\right) &\times 
\left(\varphi\circ E_D^N\circ E^M_N\right)^{it} \quad \text{or} \\
\left(\text{a reduced word in $M$}\right) &\times 
\left(\varphi\circ\theta^{-1}\circ E_{\theta(D)}^N\circ E^M_N\right)^{it}, 
\end{align*}
}and the desired assertion follows from the condition (M) for the original $M$. 
\end{proof} 

\section{Factoriality and Central Sequences} 

Let us begin by fixing our setting and notations throughout this section. Let $N$ be a $\sigma$-finite von Neumann algebra, and $\theta : D \rightarrow N$ be a normal $*$-isomorphism from a von Neumann subalgebra $D$ of $N$ into $N$. Suppose further that there are faithful normal conditional expectations 
$E_D^N : N \rightarrow D$ and $E_{\theta(D)}^N : N \rightarrow \theta(D)$. Let us consider the HNN extension: 
\begin{equation*}
\left(M, E^M_{\pi(N)}\right) = 
\left(N, E_D^N\right) \underset{D}{\bigstar} \left(\theta, E_{\theta(D)}^N\right), 
\end{equation*}
and write $M = N \underset{D}{\bigstar} \theta$ for short, when no confusion arises. The discussions in what follows heavily depend upon the construction of HNN extensions given in \S3, and thus we should briefly recall the procedure to fix notations. We first begin by  constructing the free product with amalgamation: 
\begin{equation*}
\left({\mathcal N}, {\mathcal E}\right) = 
\left(N\otimes M_2\left({\mathbf C}\right), E_{\theta} : \iota_{\theta} \right)
\underset{D\otimes{\mathbf C}^2}{\bigstar} 
\left(N\otimes M_2\left({\mathbf C}\right), E_1: \iota_1 \right), 
\end{equation*}
where 
\begin{equation*}
E_{\theta} :=
\begin{bmatrix} E_D^N &  \\  & E_{\theta(D)}^N \end{bmatrix}, \quad 
E_1 :=
\begin{bmatrix} E_D^N &  \\  & E_D^N \end{bmatrix}; 
\end{equation*}
\begin{equation*}
\iota_{\theta} :=
\begin{bmatrix} \mathrm{Id}_D &  \\  & \theta \end{bmatrix}, \quad 
\iota_1 :=
\begin{bmatrix} \mathrm{Id}_D &  \\  & \mathrm{Id}_D \end{bmatrix}.
\end{equation*}
Here, we denote the canonical embedding maps of $N\otimes M_2\left({\mathbf C}\right)$ onto the 1st/2nd free components by $\lambda_{\theta}$, $\lambda_1$, respectively, and the embedding map of $D\otimes{\mathbf C}^2$ into ${\mathcal N}$ by $\lambda$. Note that $\lambda = \lambda_{\theta}\circ\iota_{\theta} = \lambda_1\circ\iota_1$. As before, we will write 
${\mathcal N}_{\theta} := \lambda_{\theta}\left(N\otimes M_2\left({\mathbf C}\right)\right)$, ${\mathcal N}_1:= \lambda_1\left(N\otimes M_2\left({\mathbf C}\right)\right)$ and 
${\mathcal N}_{\theta}^{\circ} := {\mathcal N}_{\theta}\cap\mathrm{Ker}{\mathcal E} = \lambda_{\theta}\left(\mathrm{Ker}E_{\theta}\right)$, ${\mathcal N}_1^{\circ} := {\mathcal N}_1\cap\mathrm{Ker}{\mathcal E} = \lambda_1\left(\mathrm{Ker}E_1\right)$, and moreover denote ${\mathcal D} := \lambda\left(D\otimes{\mathbf C}^2\right)$. Set $p :=\lambda\left(1\otimes e_{11}\right) \in \lambda\left(D\otimes{\mathbf C}^2\right)$ ( $\subseteq {\mathcal N}$), and then the HNN extension $M = N \underset{D}{\bigstar} \theta$ is obtained as $M :=\langle \pi\left(N\right), u(\theta) \rangle'' \subseteq p{\mathcal N}p$ with
{\allowdisplaybreaks 
\begin{gather*}
\pi\left(n\right) :=\lambda_{\theta}\left(n\otimes e_{11}\right) = 
p\lambda_{\theta}\left(n\otimes 1\right)p, \quad n \in N;  \\
u(\theta) :=
\lambda_1\left(1\otimes e_{12}\right)\lambda_{\theta}\left(1\otimes e_{21}\right). 
\end{gather*}
}Let ${\mathcal E}_{\theta} : {\mathcal M} \rightarrow {\mathcal N}_{\theta}$ be the conditional expectation onto the 1st free component, conditioned by ${\mathcal E}$, i.e., ${\mathcal E}\circ{\mathcal E}_{\theta} = {\mathcal E}$. The conditional expectation $E_{\pi(N)}^{M} : M \rightarrow \pi(N)$ is given as the restriction of ${\mathcal E}_{\theta}$ to $M$. In this section, any embedding map appearing in the above construction will be not omitted to avoid any confusion as long as when we will treat the amalgamated free product ${\mathcal N} = {\mathcal N}_{\theta} \underset{{\mathcal D}}{\bigstar} {\mathcal N}_1$ to get any result on the HNN extension $M = N \underset{D}{\bigstar} \theta$. 

\medskip
Next, we briefly summarize some of the basics on ultraproducts of von Neumann algebras needed in this section. We refer to \cite[Chap.~5]{ocneanu:LNM} for the topic. (Also see \cite[\S\S2.2]{ueda:transAMS} as a brief summary fitting into our treatment.) Fix a free ultrafilter $\omega \in \beta({\mathbb N})\setminus{\mathbb N}$. For a given $\sigma$-finite von Neumann algebra $P$, let us denote by ${\mathcal I}_{\omega}^P$ the set of bounded sequences $\left(x_n\right)_{n\in{\mathbb N}}$ in $P$ satisfying $\sigma\text{-}s^*\text{-}\lim_{n\rightarrow\omega} x_n = 0$. With letting ${\mathcal M}\left({\mathcal I}_{\omega}^P\right) :=$ the multiplier algebra of ${\mathcal I}_{\omega}^P$ inside the algebra $\ell^{\infty}\left({\mathbb N}, P\right) = P\otimes\ell^{\infty}\left({\mathbb N}\right)$ of all bounded sequences in $P$, the ultraproduct $P^{\omega}$ is defined as the quotient $C^*$-algebra ${\mathcal M}\left({\mathcal I}_{\omega}^P\right)/{\mathcal I}_{\omega}^P$ with quotient map $\pi_{\omega}^{P}$. In what follows, we will need the following standard facts: (1) Every constant sequence in $P$ belongs to the multiplier algebra ${\mathcal M}\left({\mathcal I}_{\omega}^P\right)$. In particular, this implies that $P$ can be embedded into $P^{\omega}$ via $\pi_{\omega}^P$. (2) If $Q$ is a von Neumann subalgebra of $P$ that is the range of a faithful normal conditional expectation $E$, then the ultraproduct $Q^{\omega}$ is naturally embedded into the bigger one $P^{\omega}$ and $E$ is lifted to a faithful normal conditional expectation $E^{\omega} : P^{\omega} \rightarrow Q^{\omega}$ in the natural way. (3) For a projection $p \in P$, the reduced von Neumann algebra $p\left(P^{\omega}\right)p$ (with $p \in P \hookrightarrow P^{\omega}$) is naturally identified with $\left(pPp\right)^{\omega}$. Although it sounds trivial, one needs to care with (only) the case of infinite von Neumann algebras because of the definition of ultraproducts.   

\begin{proposition}\label{Prop5-1}
{\rm (c.f.~\cite[Lemma 2.1]{popa:adv} and \cite[Proposition 5]{ueda:transAMS})} Let ${\mathcal N} = {\mathcal N}_{\theta} \underset{{\mathcal D}}{\bigstar} {\mathcal N}_1$ be as above. Suppose that there are faithful normal states $\varphi$, $\varphi_{\theta}$ on $D$ and unitaries $v \in N_{\varphi\circ E_D^N}$, $v_{\theta} \in N_{\varphi_{\theta}\circ\theta^{-1}\circ E_{\theta(D)}^N}$ such that 
\begin{equation*}
E_D^N\left(v^n\right) = E_{\theta(D)}^N\left(v^n\right) = 0, \quad  E_{\theta(D)}^N\left(v_{\theta}^n\right) = 0
\end{equation*}
as long as $n \neq 0$. Define the state $\psi$ on $D\otimes{\mathbf C}^2$ and the unitary $V \in {\mathcal N}_{\psi\circ\lambda^{-1}\circ{\mathcal E}}$ by 
\begin{align*}
\psi\left(\mathrm{diag}\left[d_{11}, d_{22}\right]\right) 
&:=
\frac{1}{2} \left(\varphi\left(d_{11}\right) + \varphi_{\theta}\left(d_{22}\right)\right), \\
V :=
\lambda_{\theta}\left(\begin{bmatrix} v &  \\  & v_{\theta} \end{bmatrix}\right) 
&\in \lambda_{\theta}\left(
\left(N\otimes M_2\left({\mathbf C}\right)\right)_{\psi\circ\iota_{\theta}^{-1}\circ E_{\theta}}\right) 
\subseteq {\mathcal N}_{\psi\circ\lambda^{-1}\circ{\mathcal E}},  
\end{align*}
respectively. Then, we have, for every $X \in \left\langle V \right\rangle' \cap {\mathcal N}^{\omega}$, 
\begin{equation*}
\left\Vert u(\theta)\left(X - {\mathcal E}_{\theta}^{\omega}\left(X\right)\right)\right\Vert_
{L^2\left({\mathcal N}^{\omega}\right)} \leq  
\left\Vert \left[ u(\theta), X\right]\right\Vert_{L^2\left({\mathcal N}^{\omega}\right)}, 
\end{equation*}
where the canonical injection $\Lambda_{\left(\psi\circ\lambda^{-1}\circ{\mathcal E}\right)^{\omega}} : {\mathcal N}^{\omega} \rightarrow L^2\left({\mathcal N}^{\omega}\right)$ with respect to $\left(\psi\circ\lambda^{-1}\circ{\mathcal E}\right)^{\omega}$ is omitted.  
\end{proposition} 

\medskip
The idea of the proof given below is essentially the same as that of \cite[Proposition 5]{ueda:transAMS}, but not exactly the same because $u(\theta)$ is not in a single free component and indeed is in ${\mathcal N}_1^{\circ}{\mathcal N}_{\theta}^{\circ}$, a set of reduced words of length $2$. Here is a good place to mention the following: In the statement of \cite[Proposition 5]{ueda:transAMS}, it is commented that ``$uDu^* = D = wDw^*$," one of the assumptions there,  is automatic from the other one given there. This is a wrong comment, but we would like to emphasize that all the cases treated in \cite{ueda:transAMS} satisfy the condition, and moreover that the condition is never used in the proof there.

\begin{proof} Let us begin by introducing the following decomposition: 
$$
{\mathcal N}_{\theta}^{\circ} = 
{\mathcal N}_{\theta}^{\vartriangle} + 
{\mathcal N}_{\theta}^{\triangledown}, 
$$
where 
$$
{\mathcal N}_{\theta}^{\vartriangle} :=
\lambda_{\theta}\left(\begin{bmatrix} 0 & \theta(D) \\ 0 & 0 \end{bmatrix}\right), \quad 
{\mathcal N}_{\theta}^{\triangledown} :=
\lambda_{\theta}\left(\begin{bmatrix} \mathrm{Ker}E_D^N & 
\mathrm{Ker}E_{\theta(D)}^N \\ N & \mathrm{Ker}E_{\theta(D)}^N \end{bmatrix}\right).  
$$ 
Note that for each pair  
$$
X^{\vartriangle} = 
\begin{bmatrix} 0 & \theta(x) \\ 0 & 0 \end{bmatrix} \in \begin{bmatrix} 0 & \theta(D) \\ 0 & 0 \end{bmatrix}, 
\quad 
Y^{\triangledown} = 
\begin{bmatrix} y_{11} & y_{12} \\ y_{21} & y_{22} \end{bmatrix} \in 
\begin{bmatrix} \mathrm{Ker}E_D^N & 
\mathrm{Ker}E_{\theta(D)}^N \\ N & \mathrm{Ker}E_{\theta(D)}^N \end{bmatrix}, 
$$
we have 
\begin{equation}\label{Eq8} 
E_{\theta}\left(X^{\vartriangle}{}^* Y^{\triangledown}\right) 
= 
E_{\theta}\left(\begin{bmatrix} 0 & 0 \\\theta( x)^* y_{11} & \theta(x)^* y_{12} \end{bmatrix}\right) 
= 
\begin{bmatrix} 0 & 0 \\ 0 & \theta(x)^* E_{\theta(D)}^N\left(y_{12}\right) \end{bmatrix} = O.
\end{equation}
This decomposition is essential in what follows. 

In the standard Hilbert space $L^2\left({\mathcal N}\right)$, we introduce the following five subspaces: 
{\allowdisplaybreaks 
\begin{align*} 
{\mathcal X}_1 &:=
\text{the closed subspace generated by}\ 
\Lambda_{\psi\circ\lambda^{-1}\circ{\mathcal E}}\left(
{\mathcal N}_{\theta}^{\circ}\cdots{\mathcal N}_1^{\circ}\right); \\
{\mathcal X}_2 &:=
\text{the closed subspace generated by}\ 
\Lambda_{\psi\circ\lambda^{-1}\circ{\mathcal E}}\left(
{\mathcal N}_1^{\circ}\cdots{\mathcal N}_{\theta}^{\circ}\right); \\ 
{\mathcal X}_3 &:=
\text{the closed subspace generated by}\ 
\Lambda_{\psi\circ\lambda^{-1}\circ{\mathcal E}}\left(
{\mathcal N}_1^{\circ}\cdots{\mathcal N}_1^{\circ}\right); \\
{\mathcal X}_4 &:=
\text{the closed subspace generated by}\ 
\Lambda_{\psi\circ\lambda^{-1}\circ{\mathcal E}}\left(
{\mathcal N}_{\theta}^{\vartriangle}
{\mathcal N}_1^{\circ}\cdots
{\mathcal N}_{\theta}^{\circ}\right); \\
{\mathcal X}_5 &:=
\text{the closed subspace generated by}\ 
\Lambda_{\psi\circ\lambda^{-1}\circ{\mathcal E}}\left(
{\mathcal N}_{\theta}^{\triangledown}
{\mathcal N}_1^{\circ}\cdots
{\mathcal N}_{\theta}^{\circ}\right), 
\end{align*}
}where $\Lambda_{\psi\circ\lambda^{-1}\circ{\mathcal E}}$ denotes the canonical injection of ${\mathcal N}$ into $L^2\left({\mathcal N}\right)$ with respect to $\psi\circ\lambda^{-1}\circ{\mathcal E}$. Then, we have 
\begin{equation*}
L^2\left({\mathcal N}\right) = 
\left[{\mathcal X}_1\oplus{\mathcal X}_2\oplus{\mathcal X}_3\oplus
{\mathcal X}_4\right]\oplus{\mathcal X}_5\oplus L^2\left({\mathcal N}_{\theta}\right), \quad L^2\left({\mathcal N}_{\theta}\right) \overset{\mathrm{identify}}{=} \overline{\Lambda_{\psi\circ\lambda^{-1}\circ{\mathcal E}}\left({\mathcal N}_{\theta}\right)}.
\end{equation*}
Note here that ${\mathcal X}_4$ and ${\mathcal X}_5$ are orthogonal, which follows from \eqref{Eq8}. We will treat the subspaces ${\mathcal X}_1$, ${\mathcal X}_2$, ${\mathcal X}_3$, ${\mathcal X}_4$ in common, while will do ${\mathcal X}_5$ carefully by looking at $u(\theta)$.  

We introduce the operator $T_{V^n}$, $n \in {\mathbb Z}$, on 
$L^2\left({\mathcal N}\right)$ 
defined by 
\begin{equation*}
T_{V^n}
\Lambda_{\psi\circ\lambda^{-1}\circ{\mathcal E}}\left(X\right) := \Lambda_{\psi\circ\lambda^{-1}\circ{\mathcal E}}\left(V^n X V^{-n}\right), 
\quad X \in {\mathcal N}.
\end{equation*}
Since $V$ is in the centralizer of $\psi\circ\lambda^{-1}\circ{\mathcal E}$ due to Proposition \ref{Prop2-3}, one can easily 
verify: $T_{V^n}$ is a unitary; and $T_{V}^n = T_{V^n}$, $T_V^n P_{{\mathcal X}_i} = P_{T_V^n{\mathcal X}_i} T_V^n$ for every $n \in {\mathbb Z}$. Here, $P_{\mathcal Y}$ denotes the projection onto a closed subspace ${\mathcal Y}$.   

\medskip\noindent
{\bf Claim:} For $i = 1,2,3,4$ ($\neq 5$), we have 
$T_V^n{\mathcal X}_i \perp T_V^m{\mathcal X}_i$ as long as $n \neq m$. 
\begin{proof}[Proof of Claim] The first three subspaces ${\mathcal X}_1$, ${\mathcal X}_2$, ${\mathcal X}_3$ are treated exactly in the same way as in the proof of \cite[Proposition 5]{ueda:transAMS}. But we would like to explain below the case of ${\mathcal X}_1$ because it is somewhat non-trivial, and then the case of ${\mathcal X}_4$. 

Let us choose two alternating words in ${\mathcal N}_{\theta}^{\circ}$, 
${\mathcal N}_1^{\circ}$ starting at ${\mathcal N}_{\theta}^{\circ}$ 
and ending at ${\mathcal N}_1^{\circ}$: 
\begin{equation*}
W(1) = X(1)_1 Y(1)_1 \cdots Y(1)_{\ell_1}, \quad 
W(2) = X(2)_1 Y(2)_1 \cdots Y(2)_{\ell_2}
\end{equation*}
with $X(k)_j \in {\mathcal N}_{\theta}^{\circ}$ and 
$Y(k)_j \in {\mathcal N}_1^{\circ}$. Since 
\begin{equation*} 
X \in {\mathcal N}_{\theta} \Longrightarrow X - {\mathcal E}(X) \in {\mathcal N}_{\theta}^{\circ}; \quad 
Y \in {\mathcal N}_1 \Longrightarrow Y - {\mathcal E}(Y) \in {\mathcal N}_1^{\circ}
\end{equation*}
(thanks to ${\mathcal E}|_{{\mathcal N}_{\theta}} = \lambda_{\theta}\circ E_{\theta}\circ\lambda_{\theta}^{-1}$ and ${\mathcal E}|_{{\mathcal N}_1} = \lambda_1\circ E_{\theta}\circ\lambda_1^{-1}$), we have 
{\allowdisplaybreaks 
\begin{align*}
\big(
&T_V^n\Lambda_{\psi\circ\lambda^{-1}\circ{\mathcal E}}\left((W(1)\right)
\big|
T_V^m\Lambda_{\psi\circ\lambda^{-1}\circ{\mathcal E}}\left((W(2)\right)
\big)_{L^2\left({\mathcal N}\right)} \\
&= 
\psi\circ\lambda^{-1}\circ{\mathcal E}\left( 
Y(2)_{\ell_2}^* \cdots Y(2)_1^* X(2)_1^* V^{n-m} 
X(1)_1 Y(1)_1 \cdots Y(1)_{\ell_1} V^{m-n}\right) \\
&= 
\psi\circ\lambda^{-1}\circ{\mathcal E}\left( 
Y(2)_{\ell_2}^* \cdots Y(2)_1^* {\mathcal E}\left(X(2)_1^* V^{n-m} X(1)_1\right)  
Y(1)_1 \cdots Y(1)_{\ell_1} V^{m-n}\right),
\end{align*}
}by using the freeness. Iterating this procedure we finally see that 
$$
\left(T_V^n\Lambda_{\psi\circ\lambda^{-1}\circ{\mathcal E}}\left((W(1)\right)
\big|T_V^m\Lambda_{\psi\circ\lambda^{-1}\circ{\mathcal E}}\left((W(2)\right)
\right)_{L^2\left({\mathcal N}\right)} = 0
$$ 
if $\ell_1 \neq \ell_2$; or otherwise (i.e., when $\ell_1 = \ell_2 \overset{\mathrm{denote}}{=:} \ell$), we have 
{\allowdisplaybreaks 
\begin{align*}
&\left(
T_V^n\Lambda_{\psi\circ\lambda^{-1}\circ{\mathcal E}}\left((W(1)\right)
\big|
T_V^m\Lambda_{\psi\circ\lambda^{-1}\circ{\mathcal E}}\left((W(2)\right)
\right)_{L^2\left({\mathcal N}\right)} = \\
&\psi\circ\lambda^{-1}\left({\mathcal E}\left(Y(2)_{\ell}^* {\mathcal E}\left(\cdots {\mathcal E}\left(Y(2)_1^* {\mathcal E}\left(X(2)_1^* V^{n-m} X(1)_1\right)  
Y(1)_1\right) \cdots \right)Y(1)_{\ell}\right) {\mathcal E}\left(V^{m-n}\right)\right), 
\end{align*}
}and this becomes $0$ as long as $V^{m-n} \in {\mathcal N}_{\theta}^{\circ}$, i.e., $m \neq n$. Hence we are done in the case of ${\mathcal X}_1$. 

To treat the case of ${\mathcal X}_4$, it suffices to note the following 
simple fact: For each pair 
\begin{equation*}
X(1) = 
\lambda_{\theta}\left(\begin{bmatrix} 0 & \theta\left(x(1)\right) \\ 0 & 0 \end{bmatrix}\right), \quad 
X(2) = 
\lambda_{\theta}\left(\begin{bmatrix} 0 & \theta\left(x(2)\right) \\ 0 & 0 \end{bmatrix}\right)
\end{equation*}
{\allowdisplaybreaks in ${\mathcal N}_{\theta}^{\vartriangle}$, we have 
\begin{align*}
X(2)^* V^k X(1) 
&= \lambda_{\theta}\left(\begin{bmatrix} 0 & 0 \\ \theta\left(x(2)\right)^* & 0 \end{bmatrix} 
\begin{bmatrix} v^k & 0 \\ 0 & v_{\theta}^k \end{bmatrix}
\begin{bmatrix} 0 & \theta\left(x(1)\right) \\ 0 & 0 \end{bmatrix}\right) \\
&= 
\lambda_{\theta}\left(
\begin{bmatrix} 0 & 0 \\ 0 & \theta\left(x(2)\right)^* v^k \theta\left(x(1)\right) \end{bmatrix}\right) 
\in {\mathcal N}_{\theta}^{\circ}
\end{align*}
}as long as $k \neq 0$. Hence, for each pair of alternating words 
\begin{equation*}
W(1) = X(1)_1 Y(1)_1 \cdots Y(1)_{\ell_1} X(1)_{\ell_1}, \quad 
W(2) = X(2)_1 Y(2)_1 \cdots Y(2)_{\ell_2} X(2)_{\ell_2},  
\end{equation*}
with $X(k)_j \in {\mathcal N}_{\theta}^{\circ}$, $Y(k)_j \in {\mathcal N}_1^{\circ}$ ($2 \leq j \leq \ell_1$ or $\ell_2$) and $X(1)_1, X(2)_1 \in {\mathcal N}_{\theta}^{\vartriangle}$, we have 
\begin{equation*}
\big(
T_V^n\Lambda_{\psi\circ\lambda^{-1}\circ{\mathcal E}}\left((W(1)\right)
\big|
T_V^m\Lambda_{\psi\circ\lambda^{-1}\circ{\mathcal E}}\left((W(2)\right)
\big)_{L^2\left({\mathcal N}\right)} = 0
\end{equation*}
since $X(2)_1^* V^{n-m} X(1)_1 \in {\mathcal N}_{\theta}^{\circ}$ as long as 
$n \neq m$.  
\end{proof}

Let us choose an element $X = \pi_{\omega}^{\mathcal N}\left(\left(X_k\right)_{k\in{\mathbb N}}\right)$ satisfying $X = VXV^*$ with identifying $V = \pi_{\omega}^{\mathcal N}\left(\left(V,V,\dots\right)\right)$. Set 
\begin{equation*}
{\mathcal X} :=
{\mathcal X}_1\oplus{\mathcal X}_2\oplus{\mathcal X}_3\oplus{\mathcal X}_4, 
\end{equation*}
and, in the same way as in \cite[Proposition 5]{ueda:transAMS} based on the above Claim, we see that, for each $\varepsilon > 0$, there is a neighborhood $W_{\varepsilon}$ at $\omega$ (in the $w^*$-topology on the Stone-\v{C}ech compactification $\beta\left({\mathbb N}\right)$) such that 
\begin{equation*}
\left\Vert P_{\mathcal X} 
\Lambda_{\psi\circ\lambda^{-1}\circ{\mathcal E}}\left(X_k\right)
\right\Vert_{L^2\left({\mathcal N}\right)} < \varepsilon
\end{equation*}
as long as $k \in W_{\varepsilon} \cap {\mathbb N}$. For the sake of 
completeness, we will repeat the detailed argument. In what follows, we will denote $\Lambda := \Lambda_{\psi\circ\lambda^{-1}\circ{\mathcal E}}$, $\Lambda^{\omega} := \Lambda_{\left(\psi\circ\lambda^{-1}\circ{\mathcal E}\right)^{\omega}}$ for simplicity. For each fixed $n \in {\mathbb Z}$, we have 
\begin{equation*}
\lim_{k\rightarrow\omega} 
\Vert \Lambda\left(X_k - V^n X_k V^{-n}\right) \Vert_{L^2\left({\mathcal N}\right)} = 
\Vert \Lambda^{\omega}\left(X - V^n X V^{-n}\right) \Vert_{L^2\left({\mathcal N}^{\omega}\right)} 
= 0, 
\end{equation*}
and hence, for each $\delta>0$ and for each $n_0 \in {\mathbb N}$, there is a neighborhood $W$ at $\omega$ such that 
\begin{equation*}
\Vert \Lambda\left(X_k - V^n X_k V^{-n}\right) \Vert_{L^2\left({\mathcal N}\right)} < \delta 
\end{equation*}
for every $k \in W\cap{\mathbb N}$ and $n \in {\mathbb Z}$ 
with $|n| \leq n_0$, and thus for each $i \neq 5$, 
{\allowdisplaybreaks 
\begin{align*} 
\left\Vert P_{{\mathcal X}_i}\Lambda\left(X_k\right)\right\Vert_{L^2\left({\mathcal N}\right)}^2 
&= 
\left\Vert T_{V^n} P_{{\mathcal X}_i}\Lambda\left(X_k\right)\right\Vert_{L^2\left({\mathcal N}\right)}^2 \\   
&= 
\left\Vert T_{V^n} P_{{\mathcal X}_i}
\Lambda\left(X_k\right) 
- 
P_{T_{V^n}{\mathcal X}_i}
\Lambda\left(X_k\right)
+ 
P_{T_{V^n}{\mathcal X}_i}
\Lambda\left(X_k\right)
\right\Vert_{L^2\left({\mathcal N}\right)}^2 \\
&\leq 
2\left\{\left\Vert P_{T_{V^n} {\mathcal X}_i} 
\Lambda\left(V^nX_k V^{-n} - X_k\right)\right\Vert_{L^2\left({\mathcal N}\right)}^2 
+ 
\left\Vert P_{T_{V^n}{\mathcal X}_i}\Lambda\left(X_k\right)\right\Vert_{L^2\left({\mathcal N}\right)}^2\right\} \\
&< 
2\left\{\delta^2 + 
\left\Vert P_{T_{V^n}{\mathcal X}_i}\Lambda\left(X_k\right)\right\Vert_{L^2\left({\mathcal N}\right)}^2\right\}, 
\end{align*}
}and hence 
{\allowdisplaybreaks 
\begin{align*} 
\left(2n_0+1\right)\left\Vert P_{{\mathcal X}_i} 
\Lambda\left(X_k\right)
\right\Vert_{L^2\left({\mathcal N}\right)}^2 
&< 
2\left\{\left(2n_0+1\right)\delta^2 + 
\sum_{\left|n\right| \leq n_0} 
\left\Vert 
P_{T_{V^n}{\mathcal X}_i}
\Lambda\left(X_k\right)
\right\Vert_{L^2\left({\mathcal N}\right)}^2\right\} \\
&\leq 
2\left\{\left(2n_0+1\right)\delta^2 + 
\left\Vert 
\Lambda\left(X_k\right)
\right\Vert_{L^2\left({\mathcal N}\right)}^2\right\} 
\end{align*} 
}by the previous Claim. Thus, we have 
\begin{equation*}
\left\Vert P_{{\mathcal X}_i} \Lambda\left(X_k\right)\right\Vert_{L^2\left({\mathcal N}\right)}^2 \leq 
2\left\{\delta^2 + 
\frac{1}{2n_0+1}\left\Vert X_k \right\Vert_{\infty}^2\right\}
\end{equation*}
as long as $k \in W\cap{\mathbb N}$. Therefore, we get the desired assertion since ${\mathcal X}$ is the direct sum ${\mathcal X}_1\oplus {\mathcal X}_2\oplus{\mathcal X}_3\oplus{\mathcal X}_4$.  

For a while, the right $u(\theta)$ in the quantity $u(\theta)X - Xu(\theta)$ ($=[u(\theta), X]$) is replaced by an analytic element $y \in {\mathcal N}_1^{\circ}{\mathcal N}_{\theta}^{\circ}$ under the modular action $\sigma_t^{\psi\circ\lambda^{-1}\circ{\mathcal E}}$. (Note that the restriction of $\sigma_t^{\psi\circ\lambda^{-1}\circ{\mathcal E}}$ to ${\mathcal N}_1^{\circ}{\mathcal N}_{\theta}^{\circ}$ is nothing but the ``product $\sigma_t^{\varphi\circ\iota_1\circ E_1}(\ \cdot\ ) \sigma_t^{\varphi\circ\iota_{\theta}\circ E_{\theta}}(\ \cdot\ )$.'') 

The Hilbert space $L^2\left({\mathcal N}^{\omega}\right)$ can be isometrically embedded into the ultraproduct Hilbert space $L^2\left({\mathcal N}\right)^{\omega}$ (see e.g.~\cite[\S\S2.2]{ueda:transAMS}), and the embedding is given by 
$\Lambda^{\omega}\left(
\pi_{\omega}^{\mathcal N}\left(
\left(x_n\right)_{n\in{\mathbb N}}
\right)
\right) 
\in L^2\left({\mathcal N}^{\omega}\right) 
\mapsto 
\left[
\left(
\Lambda
\left(x_n\right)
\right)_{n\in{\mathbb N}}
\right]_{L^2\left({\mathcal N}\right)^{\omega}} 
\in L^2\left({\mathcal N}\right)^{\omega}$, the quotient class of the given bounded sequence $\left(\Lambda\left(x_n\right)\right)_{n\in{\mathbb N}} \in \ell^{\infty}\left({\mathbb N}, L^2\left({\mathcal N}\right)\right)$. Hence, we will regard $L^2\left({\mathcal N}^{\omega}\right)$ as a closed subspace of $L^2\left({\mathcal N}\right)^{\omega}$ via the embedding. 

We have 
{\allowdisplaybreaks 
\begin{align*} 
&\left\Vert 
\Lambda^{\omega}
\left(
u(\theta)\left(X - {\mathcal E}_{\theta}^{\omega}\left(X\right)\right)
\right) - 
\left[\left(u(\theta) P_{{\mathcal X}_5}
\Lambda
\left(X_k\right)\right)_{k\in{\mathbb N}}\right]_{L^2\left({\mathcal N}\right)^{\omega}} 
\right\Vert^2_{L^2\left({\mathcal N}\right)^{\omega}} \\
&\leq 
\sup_{k \in W_{\varepsilon}\cap{\mathbb N}} 
\left\Vert u(\theta)\left(
1_{L^2\left({\mathcal N}\right)} 
- 
P_{L^2\left({\mathcal N}_{\theta}\right)} -P_{{\mathcal X}_5}
\right)
\Lambda
\left(X_k\right) 
\right\Vert_{L^2\left({\mathcal N}\right)}^2 \\
&\leq 
\sup_{k \in W_{\varepsilon}\cap{\mathbb N}} 
\left\Vert P_{\mathcal X}
\left(X_k\right) 
\right\Vert_{L^2\left({\mathcal N}\right)}^2 \leq \varepsilon,    
\end{align*}
}and hence 
\begin{equation*}
\Lambda^{\omega}
\left(
u(\theta)\left(X - {\mathcal E}_{\theta}^{\omega}\left(X\right)\right)
\right) 
= 
\left[\left(u(\theta) P_{{\mathcal X}_5}
\Lambda
\left(X_k\right)\right)_{k\in{\mathbb N}}\right]_{L^2\left({\mathcal N}\right)^{\omega}}  
\end{equation*}
inside $L^2\left({\mathcal N}\right)^{\omega}$ since $\varepsilon$ is arbitrary. With the notations $\displaystyle{\sigma_t := \sigma_t^{\psi\circ\lambda^{-1}\circ{\mathcal E}}}$, $\displaystyle{J := J_{\varphi\circ\lambda^{-1}\circ{\mathcal E}}}$, we compute 
{\allowdisplaybreaks 
\begin{align*} 
&\left\Vert 
\Lambda^{\omega}
\left(
\left(X - {\mathcal E}_{\theta}^{\omega}\left(X\right)\right)y
\right) - 
\left[\left(
J\sigma_{-\frac{\sqrt{-1}}{2}}\left(y^*\right)J
P_{{\mathcal X}_5}
\Lambda
\left(X_k\right)\right)_{k\in{\mathbb N}}\right]_{L^2\left({\mathcal N}\right)^{\omega}} 
\right\Vert^2_{L^2\left({\mathcal N}\right)^{\omega}} \\
&\leq 
\sup_{k \in W_{\varepsilon}\cap{\mathbb N}} 
\left\Vert 
J\sigma_{-\frac{\sqrt{-1}}{2}}\left(y^*\right)J 
\Lambda\left(
X_k - {\mathcal E}_{\theta}\left(X_k\right)
\right) 
- 
J\sigma_{-\frac{\sqrt{-1}}{2}}\left(y^*\right)J 
P_{{\mathcal X}_5}
\Lambda
\left(X_k\right) 
\right\Vert_{L^2\left({\mathcal N}\right)}^2 \\ 
&\leq \left\Vert\sigma_{-\frac{\sqrt{-1}}{2}}\left(y^*\right)\right\Vert_{\infty}\cdot\varepsilon.     
\end{align*}
}Since $\varepsilon$ is arbitrary, we get 
\begin{equation*}
\Lambda^{\omega}
\left(
\left(X - {\mathcal E}_{\theta}^{\omega}\left(X\right)\right)y
\right) = 
\left[\left(
J\sigma_{-\frac{\sqrt{-1}}{2}}\left(y^*\right)J
P_{{\mathcal X}_5}
\Lambda
\left(X_k\right)\right)_{k\in{\mathbb N}}\right]_{L^2\left({\mathcal N}\right)^{\omega}}
\end{equation*}
inside $L^2\left({\mathcal N}\right)^{\omega}$. 
We also have
\begin{equation*}
\Lambda^{\omega}
\left(
u(\theta)
\left(X - {\mathcal E}_{\theta}^{\omega}\left(X\right)\right)y
\right) = 
\left[\left(
\Lambda
\left(
u(\theta) {\mathcal E}_{\theta}\left(X_k\right) - 
{\mathcal E}_{\theta}\left(X_k\right)y
\right)
\right)_{k\in{\mathbb N}}\right]_{L^2\left({\mathcal N}\right)^{\omega}} 
\end{equation*}
inside $L^2\left({\mathcal N}\right)^{\omega}$. Therefore, it suffices to show that 
\begin{equation*}
u(\theta)\left({\mathcal N}_{\theta}^{\triangledown}
{\mathcal N}_1^{\circ}\cdots
{\mathcal N}_{\theta}^{\circ}\right), \quad 
\left({\mathcal N}_{\theta}^{\triangledown}
{\mathcal N}_1^{\circ}\cdots
{\mathcal N}_{\theta}^{\circ}\right)\left({\mathcal N}_1^{\circ}{\mathcal N}_{\theta}^{\circ}\right), \quad 
u(\theta){\mathcal N}_{\theta} + 
{\mathcal N}_{\theta}\left({\mathcal N}_1^{\circ}{\mathcal N}_{\theta}^{\circ}\right) 
\end{equation*}
are mutually orthogonal with respect to $\psi\circ\lambda^{-1}\circ{\mathcal E}$. Indeed, if this assertion was true, then it would 
follow that  
\begin{equation*}
u(\theta) P_{{\mathcal X}_5}\Lambda\left(X_k\right), \quad 
J\sigma_{-\frac{\sqrt{-1}}{2}}\left(y^*\right)J
P_{{\mathcal X}_5}\Lambda\left(X_k\right), \quad
\Lambda\left(u(\theta) {\mathcal E}_{\theta}\left(X_k\right) - 
{\mathcal E}_{\theta}\left(X_k\right)y\right)
\end{equation*}
are mutually orthogonal in $L^2\left({\mathcal N}\right)$ for every $k$, and so are 
\begin{equation*}
\Lambda^{\omega}
\left(u(\theta)
\left(X-{\mathcal E}_{\theta}^{\omega}\left(X\right)\right) 
\right), \quad 
\Lambda^{\omega}\left(
\left({\mathcal E}_{\theta}^{\omega}\left(X\right) 
- X\right)u(\theta)\right), \quad
\Lambda^{\omega}
\left(u(\theta){\mathcal E}_{\theta}^{\omega}\left(X\right) - 
{\mathcal E}_{\theta}^{\omega}\left(X\right)u(\theta) \right)  
\end{equation*} 
in $L^2\left({\mathcal N}^{\omega}\right)$ (or more precisely, inside $L\left({\mathcal N}\right)^{\omega}$) since one can find a bounded net of analytic elements in ${\mathcal N}_1^{\circ}{\mathcal N}_{\theta}^{\circ}$ that converges to $u(\theta)$ in the $\sigma$-strong$^*$ topology.  

Let us now show the desired orthogonal relation among the words 
in question. For an alternating word 
\begin{equation*} 
W = X_1 Y_2 \cdots X_{\ell} \in {\mathcal N}_{\theta}^{\triangledown}{\mathcal N}_1^{\circ}\cdots{\mathcal N}_{\theta}^{\circ}
\end{equation*}
with 
\begin{equation*}
X_1 = 
\lambda_{\theta}\left(\begin{bmatrix} x_{11} & x_{12} \\ x_{21} & x_{22} \end{bmatrix}\right) 
\in {\mathcal N}_{\theta}^{\triangledown} = 
\begin{bmatrix} \mathrm{Ker}E_D^N & \mathrm{Ker}E_{\theta(D)}^N \\ 
N & \mathrm{Ker}E_{\theta(D)}^N \end{bmatrix}, 
\end{equation*}
we have 
{\allowdisplaybreaks 
\begin{align*}
u(\theta)W &= 
\lambda_1\left(\begin{bmatrix} 0 & 1 \\ 0 & 0 \end{bmatrix}\right) 
\lambda_{\theta}\left(\begin{bmatrix} 0 & 0 \\ 1 & 0 \end{bmatrix}\right) 
\lambda_{\theta}\left(\begin{bmatrix} x_{11} & x_{12} \\ x_{21} & x_{22} \end{bmatrix}\right) 
Y_2 \cdots \\
&= 
\lambda_1\left(\begin{bmatrix} 0 & 1 \\ 0 & 0 \end{bmatrix}\right) 
\lambda_{\theta}\left(\begin{bmatrix} 0 & 0 \\ x_{11} & x_{12} \end{bmatrix}\right) Y_2 \cdots
\in {\mathcal N}_1^{\circ}{\mathcal N}_{\theta}^{\circ}{\mathcal N}_1^{\circ} \cdots,  
\end{align*}
}and hence  
\begin{equation*}
u(\theta)\big(\underbrace{{\mathcal N}_{\theta}^{\triangledown}
{\mathcal N}_1^{\circ}\cdots
{\mathcal N}_{\theta}^{\circ}}_{\text{length $\ell_1$}}\big)
 \subseteq
{\mathcal N}_1^{\circ}{\mathcal N}_{\theta}^{\circ}{\mathcal N}_1^{\circ}\cdots {\mathcal N}_{\theta}^{\circ}\  
\text{is of length $\ell_1 + 1 \geq 4$}. 
\end{equation*}
Since 
\begin{equation*}
\big(\underbrace{{\mathcal N}_{\theta}^{\triangledown}
{\mathcal N}_1^{\circ}\cdots
{\mathcal N}_{\theta}^{\circ}}_{\text{length $ \ell_2$}}\big)\left({\mathcal N}_1^{\circ}{\mathcal N}_{\theta}^{\circ}\right)\ \text{is of length $\ell_2 +2 \geq 5$}, 
\end{equation*} 
we see that the first two sets of words in question are orthogonal. Notice here that the length of every reduced word living in the set  
\begin{equation}\label{Eq9}
u(\theta){\mathcal N}_{\theta} + 
{\mathcal N}_{\theta}\left({\mathcal N}_1^{\circ}{\mathcal N}_{\theta}^{\circ}\right)
\subseteq 
{\mathcal N}_1^{\circ} + 
{\mathcal N}_1^{\circ}{\mathcal N}_{\theta}^{\circ} + 
{\mathcal N}_{\theta}^{\circ}{\mathcal N}_1^{\circ}{\mathcal N}_{\theta}^{\circ}
\end{equation}
is less than 3, and hence the left-hand side of \eqref{Eq9} is easily seen (by looking at the lengths of words) to be orthogonal to the other two sets of words in question. Hence we are done.  
\end{proof} 

Here, we should give simple facts concerning ultraproducts: (a) It is easy to see that $\pi : N \rightarrow M$ and $\theta : D \rightarrow N$ can be lifted to the normal $*$-isomorphisms $\pi^{\omega} : N^{\omega} \rightarrow M^{\omega}$ and $\theta^{\omega} : D^{\omega} \rightarrow N^{\omega}$ in the obvious manner (note that $\pi(N)$ and $\theta(D)$ are the ranges of faithful normal conditional expectations), and it follows from their construction that $\pi^{\omega}\left(N^{\omega}\right) = \pi(N)^{\omega}$ (inside $M^{\omega}$) and $\theta^{\omega}\left(D^{\omega}\right) = \theta(D)^{\omega}$ (inside $N^{\omega}$). (b) We have $\pi^{\omega}|_N = \pi$ via the embeddings $M \hookrightarrow M^{\omega}$ and $N \hookrightarrow N^{\omega}$. We will use these facts with no explicit explanation in what follows.      

\begin{theorem}\label{Thm5-2} 
Suppose that there are faithful normal states $\varphi$, $\varphi_{\theta}$ on $D$ 
and unitaries $v \in N_{\varphi\circ E_D^N}$, $v_{\theta} \in N_{\varphi_{\theta}\circ\theta^{-1}\circ E_{\theta(D)}^N}$ such that 
\begin{equation*}
E_D^N\left(v^n\right) = E_{\theta(D)}^N\left(v^n\right) = 0, \quad  E_{\theta(D)}^N\left(v_{\theta}^n\right) = 0
\end{equation*}
as long as $n \neq 0$. Then we have 
\begin{equation*}
\langle \pi\left(v\right), u(\theta) \rangle' \cap M^{\omega} 
\subseteq \pi\left(N\right)^{\omega} = \pi^{\omega}\left(N^{\omega}\right) \cong N^{\omega}\ (\text{via $\pi^{\omega}$}). 
\end{equation*}
The same relative commutant property still holds when replacing $M^{\omega}$ and $\pi(N)^{\omega}$ by $M$ and $\pi(N)$, respectively. Namely, we have
\begin{equation*}
\langle \pi\left(v\right), u(\theta) \rangle' \cap M 
\subseteq \pi\left(N\right) \cong N\ (\text{via $\pi$}). 
\end{equation*}    
\end{theorem} 
\begin{proof} Notice that 
$\pi\left(N\right) = p{\mathcal N}_{\theta}p \subseteq M \subseteq p{\mathcal N}p$, i.e., 
$\pi\left(N\right)^{\omega} = \left(p{\mathcal N}_{\theta}p\right)^{\omega} = 
p{\mathcal N}_{\theta}^{\omega}p \subseteq M^{\omega} \subseteq \left(p{\mathcal N}p\right)^{\omega} = 
p{\mathcal N}^{\omega}p$ and that $\pi\left(v\right) = pV = Vp = pVp$. Hence, the first assertion follows from Proposition \ref{Prop5-1}. 
In fact, we have 
\begin{equation*}
\langle \pi\left(v\right), u(\theta) \rangle' \cap M^{\omega} \subseteq 
\langle V, u(\theta) \rangle' \cap p{\mathcal N}^{\omega}p \subseteq 
p{\mathcal N}^{\omega}p = \pi\left(N\right)^{\omega}
\end{equation*}
since $u(\theta)^* u(\theta) = u(\theta)u(\theta)^* = p$. 

Let us choose $x \in M \cap \pi(N)^{\omega}$ inside $M^{\omega}$, and then we get  
\begin{equation*}
x = \left(E_{\pi(N)}^M\right)^{\omega}\left(x\right) = \left[\left(E_{\pi(N)}^M(x),E_{\pi(N)}^M(x),\dots\right)\right] = E_{\pi(N)}^M(x) \in \pi(N) 
\ \left(\text{inside $M^{\omega}$}\right),
\end{equation*}
Hence, the last assertion follows.    
\end{proof} 

\begin{corollary}\label{Cor5-3} 
Under the same assumption as in Theorem \ref{Thm5-2}, we have 
\begin{gather*}
{\mathcal Z}\left(M\right) = 
\langle u(\theta) \rangle' \cap \pi\left({\mathcal Z}\left(N\right)\right) \cong \left\{ x \in D\cap\theta(D)\cap N' : \theta(x) =  x \right\}\ (\text{via $\pi$}), \\
M' \cap M^{\omega} = M' \cap \pi(N)^{\omega} \cong \left\{ x \in D^{\omega} \cap\theta^{\omega}\left(D^{\omega}\right)\cap N' : \theta^{\omega}(x) =  x \right\}\ (\text{via $\pi^{\omega}$}). 
\end{gather*}
Therefore, if $N$ is further assumed to be a factor, then so is the HNN extension $M$. Moreover, the same is true for the continuous core, that is, 
\begin{equation*}
{\mathcal Z}\left(\widetilde{M}\right) = \langle u(\theta) \rangle' \cap 
\widetilde{\pi}\left({\mathcal Z}\left(\widetilde{N}\right)\right) \cong 
\left\{ x \in \widetilde{D}\cap\widetilde{\theta}\left(\widetilde{D}\right)\cap \widetilde{N}' : 
\widetilde{\theta}(x) = x \right\}\ (\text{via $\widetilde{\pi}$}), 
\end{equation*}
where $\widetilde{M}\supseteq\widetilde{N}\supseteq\widetilde{D}$, $\widetilde{\theta} : \widetilde{D} \rightarrow \widetilde{N}$, etc., are as in \S4. Thus, the flow of weights of $M$ is a factor flow of that of $N$. 
\end{corollary} 
\begin{proof} Thanks to Theorem \ref{Thm5-2}, it suffices to show the following:  
\begin{itemize} 
\item[(i)] $\langle u(\theta) \rangle' \cap \pi\left({\mathcal Z}\left(N\right)\right) \cong \left\{ x \in D\cap\theta(D)\cap N' : \theta(x) =  x \right\}$ via $\pi$. 
\item[(ii)] $\left\langle u(\theta) \right\rangle' \cap \pi(N)' \cap \pi^{\omega}\left(N^{\omega}\right) \cong \left\{ x \in D^{\omega} \cap\theta^{\omega}\left(D^{\omega}\right)\cap N' : \theta^{\omega}(x) =  x \right\}$ via $\pi^{\omega}$. 
\item[(iii)] One can construct faithful normal states $\widetilde{\varphi}$, $\widetilde{\varphi_{\theta}}$ 
on $\widetilde{D}$ in such a way that $v$ and $v_{\theta}$ are in 
the centralizers of $\widetilde{\varphi}\circ\widehat{E_D^N}$ and $\widetilde{\varphi_{\theta}}
\circ\widetilde{\theta}^{-1}\circ\widehat{E_{\theta(D)}^N}$, respectively. 
\end{itemize}
Note here that the continuous core $\widetilde{M}$ can be written again the HNN extension of the base algebra $\widetilde{N}$ by $\widetilde{\theta} : \widetilde{D} \rightarrow \widetilde{N}$ with respect to $\widehat{E_D^N}$ and $\widehat{E_{\theta(D)}^N}$ thanks to Theorem \ref{Thm4-5} and also that $\widehat{E_D^N}\left(v^n\right) = \widehat{E_{\theta(D)}^N}\left(v^n\right) = \widehat{E_{\theta(D)}^N}\left(v_{\theta}^n\right) = 0$ as long as $n \neq 0$ since $\widehat{E_D^N}\big|_N = E_D^N$ and $\widehat{E_{\theta(D)}^N}\big|_N = E_{\theta(D)}^N$. Thus, the above (iii) is indeed enough to complete the proof of the assertion on the continuous core $\widetilde{M}$. 

Let us choose $x \in {\mathcal Z}\left(N\right)$ with $u(\theta)\pi(x)u(\theta)^* = \pi(x)$ 
(and hence $u(\theta)^*\pi(x)u(\theta) = \pi(x)$). 
Then, the characterization of HNN extensions, i.e., Theorem \ref{Thm3-3} (especially, the condition (M)), enables us to compute 
{\allowdisplaybreaks 
\begin{align*}
\pi(x) 
&= E_{\pi(N)}^M\left(u(\theta)\pi(x)u(\theta)^*\right) \\
&= E_{\pi(N)}^M\left(u(\theta)\pi\left(x - E_{\theta(D)}^N(x)\right)u(\theta)^*\right) + 
E_{\pi(N)}^M\left(\pi\left(\theta^{-1}\left(E_{\theta(D)}^N(x)\right)\right)\right) \\
&= \pi\left(\theta^{-1}\left(E_{\theta(D)}^N(x)\right)\right) \in \pi(D), 
\end{align*} 
}and similarly $\pi(x) = E_{\pi(N)}^M\left(u(\theta)^*\pi(x)u(\theta)\right) = \pi\left(\theta\left(E_D^N(x)\right)\right) \in \pi\left(\theta(D)\right)$. These imply the desired assertion (i). 

The assertion (ii) is shown in the same manner, but the reader should notice the following two simple facts: (a) $u(\theta)\pi^{\omega}\left(\theta^{\omega}(x)\right)u(\theta)^* = \pi^{\omega}\left(x\right)$ for every $x \in D^{\omega}$. (b) The restriction of $\left(E_{\pi(N)}^M\right)^{\omega}$ to $*\text{-}\mathrm{Alg}\left\langle \pi^{\omega}\left(N^{\omega}\right), u(\theta) \right\rangle$ satisfies the condition (M), where $*\text{-}\mathrm{Alg}\left\langle \pi^{\omega}\left(N^{\omega}\right), u(\theta)\right\rangle$ denotes the $*$-algebra algebraically generated by $\pi^{\omega}\left(N^{\omega}\right)$ and $u(\theta)$. (Concerning (ii), we do not know whether or not $M^{\omega} \cong N^{\omega} \underset{D^{\omega}}{\bigstar} \theta^{\omega}$ since it is highly non-trivial whether or not $M^{\omega}$ is generated by $\pi^{\omega}\left(N^{\omega}\right)$ and $u(\theta)$ as von Neumann algebra. Probably, ``No!")  

We will finally prove the desired assertion (iii). Let us consider the faithful normal conditional expectations 
{\allowdisplaybreaks 
\begin{align*} 
\widehat{\varphi} &: \widetilde{D} \left(\subseteq \widetilde{N}\right) \rightarrow 
L_{\varphi}\left({\mathbf R}\right) :=
\big\langle \left(\varphi\circ E_D^N\right)^{it} (t \in {\mathbf R}) \big\rangle'' 
\cong L\left({\mathbf R}\right), \\ 
\widehat{\varphi_{\theta}} &: 
\widetilde{D} \left(\subseteq \widetilde{N}\right) \rightarrow 
L_{\varphi_{\theta}}\left({\mathbf R}\right) :=
\big\langle 
\left(\varphi_{\theta}\circ E_D^N\right)^{it} (t \in {\mathbf R}) 
\big\rangle'' \cong L\left({\mathbf R}\right)
\end{align*} 
}constructed in such a way that  
{\allowdisplaybreaks 
\begin{align*}
\widehat{\varphi}\left(\int_{-\infty}^{\infty} x(t) \left(\varphi\circ E_D^N\right)^{it} dt\right) &= 
\int_{-\infty}^{\infty} \varphi\left(x(t)\right) \left(\varphi\circ E_D^N\right)^{it} dt, \\ 
\widehat{\varphi_{\theta}}\left(\int_{-\infty}^{\infty} y(t) \left(\varphi_{\theta}\circ E_D^N\right)^{it} dt\right) &= 
\int_{-\infty}^{\infty} \varphi_{\theta}\left(y(t)\right) 
\left(\varphi_{\theta}\circ E_D^N\right)^{it} dt 
\end{align*}
}for ``smooth" functions $x(t)$, $y(t) : {\mathbf R} \rightarrow D$. We then construct two faithful normal states $\widetilde{\varphi} :=\psi\circ\widehat{\varphi}$, $\widetilde{\varphi_{\theta}} :=\psi_{\theta}\circ\widehat{\varphi_{\theta}}$ 
on $\widetilde{D}$ with faithful normal states $\psi$ and $\psi_{\theta}$ on $L_{\varphi}\left({\mathbf R}\right)$ and $L_{\varphi_{\theta}}\left({\mathbf R}\right)$, respectively. These are the desired ones. In fact, we have 
{\allowdisplaybreaks \begin{align*}
\widetilde{\varphi_{\theta}}\circ\widetilde{\theta}^{-1}\circ \widehat{E_{\theta(D)}^N}&\left(v_{\theta} \left(\int_{-\infty}^{\infty} x(t) \left(\varphi_{\theta}\circ\theta^{-1}\circ E_{\theta(D)}^N\right)^{it} dt\right)\right) \\
&= \widetilde{\varphi_{\theta}}\circ\widetilde{\theta}^{-1}\circ \widehat{E_{\theta(D)}^N}\left(\left(\int_{-\infty}^{\infty} x(t)\left(\varphi_{\theta}\circ\theta^{-1}\circ E_{\theta(D)}^N\right)^{it} dt\right) v_{\theta} \right).
\end{align*}}This follows from the assumption that $v_{\theta}$ is in the centralizer of $\varphi_{\theta}\circ \theta^{-1}\circ E_{\theta(D)}^N$. Hence, $v_{\theta}$ is in the centralizer of $\widetilde{\varphi_{\theta}}\circ\widetilde{\theta}^{-1}\circ\widehat{E^N_{\theta(D)}}$ too. The other case is quite similar and easier to show, and thus left to the reader. Hence, the desired assertion (iii) is verified. 
\end{proof}   

\begin{remark}\label{Rem5-4}{\bf [Analog of Higman, Neumann and Neumann's theorem]} {\rm Following G.~Higman, B.H.~Neumann and H.~Neumann  \cite{higmanneumann^2:JLondonMath1949} (see also  \cite[Theorem 3.1 in p.188]{lyndonschupp:book}), we can show the following:} Each finite von Neumann algebra $P$ with separable predual has a full type II$_1$ factor $\widetilde{P}$ generated by two Haar unitaries, into which $P$ can be embedded. {\rm What is new is that the generators of $\widetilde{P}$ are chosen to be Haar unitaries. In fact, for a given $P$, there are several ways based on known results to construct such a bigger type II$_1$ factor with two unitary generators, and moreover it can be made to be full or to have the Property $\Gamma$. However, the construction given below remains to work even in the $C^*$-algebra setting. We should also point out that Connes' approximate embedding problem (see \cite[lines 13-9 from the bottom in p.105]{connes:annals76}) can be read, from the viewpoint here, as whether or not any possible set of ``relations" between two Haar unitaries can be realized in the ultraproduct $R^{\omega}$ of the AFD type II$_1$ factor.} 
\end{remark}   
\begin{proof}[Proof of the assertion in Remark \ref{Rem5-4}] Since $P$ has the separable predual, it has an at most countable generating set of unitaries, say $\left\{u_n\right\}_{n \in {\mathbb Z}}$ with $u_0 = 1$. We choose two copies $S$, $T$ of the free group factor $L\left({\mathbb F}_2\right)$ with $*$-free Haar unitary generators $a,b$ and $c,d$, respectively. Choose a faithful normal tracial state $\tau_P$ on $P$, and let $\tau_S$ and $\tau_T$ be the unique tracial states on $S$ and $T$, respectively. We first embed $P$ into the free product with amalgamation:   
\begin{equation*}
N :=\left(P \bigstar S\right) 
\underset{L\left({\mathbb F}_{\infty}\right)}{\bigstar} T
\end{equation*}
with respect to $\tau_P$, $\tau_S$ and $\tau_T$ (or more precisely, the conditional expectations determined by them), where the amalgamation is taken by the identification $P \bigstar S \ni u_n\left(b^n a b^*{}^n\right) \leftrightarrow d^n c d^*{}^n \in T$. In fact, it is known that the $b^n a b^*{}^n$'s and the $d^n c d^*{}^n$'s form $*$-free families of Haar unitaries in $S$ and $T$, respectively, and furthermore it can be easily verified that the $u_n b^n a b^*{}^n$'s also form a $*$-free family of Haar unitaries in $P \bigstar S$ since the $u_n$'s and the $b^n a b^*{}^n$'s are chosen from different free components of $P \bigstar S$. Thus, the $u_n b^n a b^*{}^n$'s and the $d^n c d^*{}^n$'s generate two copies $Q_1$ and $Q_2$ of the free group factor $L\left({\mathbb F}_{\infty}\right)$ inside $P \bigstar S$ and $T$, respectively. Therefore, the above amalgamation procedure is justified and agrees with the tracial states $\tau_P\bigstar\tau_S$ and $\tau_T$, and hence the free product state $\tau_N := \left(\tau_P\bigstar\tau_S\right)\bigstar\tau_T$ becomes a trace (see \cite[\S3]{popa:invent:irredindex>4}). Since $u_0 = 1$, one has $a = c$ and $u_n = \left(d^n c d^*{}^n\right)\left(b^n a^* b^*{}^n\right) = d^n c d^*{}^n b^n c^* b^*{}^n \in \langle b, c, d \rangle''$ so that $N$ is generated by three Haar unitaries $b$, $c$, $d$. By using the normal $*$-isomorphism $\theta$ from $S = \langle b, c=a \rangle''$ onto $T = \langle c, d \rangle''$ given by $b \mapsto c$ and $c \mapsto d$, we finally embed $P$ ($\hookrightarrow N$) into the HNN extension 
\begin{equation*}
\widetilde{P} :=
N \underset{L\left({\mathbb F}_2\right)}{\bigstar} \theta 
= \langle N, u(\theta) \rangle'' 
\end{equation*}
with respect to the $\tau_N$-preserving conditional expectations. Since $c = u(\theta) d u(\theta)^*$ and $b = u(\theta)^2 d u(\theta)^*{}^2$, we see that $\widetilde{P} = \langle d, u(\theta) \rangle''$. It remains only to show that $\widetilde{P}$ is a full type II$_1$ factor. 
 
Since $S$ and $T$ are isomorphic type II$_1$ factors, we get $\tau_S = \tau_T\circ\theta$, and hence Corollary \ref{Cor4-2} says that $\tau_N$ is extended to a tracial state on $\widetilde{P}$ by the canonical conditional expectation from $\widetilde{P}$ onto $N$ since $\tau_N$ agrees with $\tau_S$, $\tau_T$. Set $v := bd$, a unitary in $N$. It is not hard to verify that $b^n$ ($n\neq0$) is orthogonal to $Q_1$ with respect to the tracial state $\tau_P\bigstar\tau_S$. It is also known that the canonical unitary generator of ${\mathbb Z}$ in the crossed-product description $T = Q_2\rtimes{\mathbb Z}$ is given by $d$ (see \cite[Proposition 4.1 and Corollary 4.2]{phillips:duke76}). These two facts show that $v = bd$ is in $\left(\left(P\bigstar S\right)\ominus Q_1\right)\left(T\ominus Q_2\right)$, a set of reduced word of length $2$ in the free product with amalgamation 
$\displaystyle{N=\left(P\bigstar S \supseteq Q_1\right)\bigstar\left(T \supseteq Q_2\right)}$. Therefore, $E_S^N\left(v^n\right) =  E_T^N\left(v^n\right) = 0$ as long as $n \neq 0$. Hence, Corollary \ref{Cor5-3} implies that $\widetilde{P}$ is a type II$_1$ factor and $\left(\widetilde{P}\right)' \cap \left(\widetilde{P}\right)^{\omega} \subseteq N'\cap N^{\omega}$. The above-mentioned two facts on $b, d$ also say that the unitaries $b, d$ satisfy the necessary conditions to apply \cite[Proposition 5]{ueda:transAMS} to $\displaystyle{N=\left(P\bigstar S \supseteq Q_1\right)\bigstar \left(T \supseteq Q_2\right)}$ ({\it n.b.}, the regularity conditions $bQ_1 b^* = Q_1$, $dQ_2 d^* =Q_2$ are not needed, see the comment given just below the statement of Proposition \ref{Prop5-1}), and therefore $N' \cap N^{\omega} \subseteq L\left({\mathbb F}_{\infty}\right)' \cap L\left({\mathbb F}_{\infty}\right)^{\omega} = {\mathbf C}1$ since $Q_1=Q_2 = L\left({\mathbb F}_{\infty}\right)$. 
\end{proof} 

\section{Several Concrete Settings} 

To illustrate how the results obtained in \S5 can be applied, we will investigate HNN extensions of von Neumann algebras in three kinds of concrete settings.  

\subsection{The HNN Extensions associated with Non-commutative Tori}\label{SS6-1} \subsubsection{Setting.} Let $\alpha \in [0,1)$ be an irrational number, and the non-commutative torus $C({\mathbb T}^2_{\alpha})$ is the universal $C^*$-algebra generated by two unitaries $u_{\alpha}$, $v_{\alpha}$ with $u_{\alpha}v_{\alpha} = e^{2\pi\sqrt{-1}\alpha}v_{\alpha}u_{\alpha}$. It is known that there is a unique tracial state $\tau_{\alpha}$ on $C({\mathbb T}_{\alpha}^2)$ determined by $\tau_{\alpha}\left(u_{\alpha}^n v_{\alpha}^m\right) = \delta_{n 0}\, \delta_{m 0}$, and its GNS representation gives the AFD type II$_1$ factor $R$, i.e., $R = C({\mathbb T}^2_{\alpha})''$ (in the representation) with the unique tracial state $\tau_R$, which is the natural extension of $\tau_{\alpha}$ to $R$. Notice that the generating unitaries $u_{\alpha}$ and $v_{\alpha}$ generate respectively two distinguished Cartan subalgebras $D_{\alpha}:=\left\langle u_{\alpha}\right\rangle''$ and $C_{\alpha} :=\left\langle v_{\alpha}\right\rangle''$ in $R$, which are both isomorphic to $L^{\infty}({\mathbb T})$. Then, let us define the normal $*$-isomorphism $\theta_{\alpha} : D_{\alpha} \rightarrow C_{\alpha} \subseteq R$ in such a way that $\theta_{\alpha}\left(u_{\alpha}\right) = v_{\alpha}$. There are two unique normal conditional expectations $E_{D_{\alpha}}^R$, $E_{\theta_{\alpha}\left(D_{\alpha}\right)}^R =E_{C_{\alpha}}^R$ from $R$ onto $D_{\alpha}$ and $\theta_{\alpha}\left(D_{\alpha}\right) = C_{\alpha}$, respectively, both of which preserve the trace $\tau_R$, and we then construct the HNN extension 
\begin{equation*}
\left(M_{\alpha}, E_R^{M_{\alpha}}\right) = 
\left(R, E_{D_{\alpha}}^R\right) \underset{D_{\alpha}}{\bigstar}\left(\theta_{\alpha}, E_{\theta_{\alpha}\left(D_{\alpha}\right)}^R\right)
\end{equation*}
with stable unitary $u(\alpha) :=u\left(\theta_{\alpha}\right)$, where the canonical embedding map of $R$ into $M$ is omitted.

\subsubsection{Trace.} It is plain to see that 
\begin{equation}\label{Eq10} 
\tau_R = \left(\tau_R\big|_{D_{\alpha}}\right)\circ E_{D_{\alpha}}^R = \left(\tau_R\big|_{D_{\alpha}}\right)\circ\theta_{\alpha}^{-1}\circ E_{\theta_{\alpha}\left(D_{\alpha}\right)}^R.  
\end{equation}
Thus, by Corollary \ref{Cor4-2}, $\tau_{M_{\alpha}} :=\tau_R \circ E_R^{M_{\alpha}}$ becomes a tracial state. 

\subsubsection{Factoriality and Fullness.} With letting $v :=u_{\alpha}v_{\alpha} = e^{2\pi\sqrt{-1}\alpha}v_{\alpha}u_{\alpha}$, it is plain to see that $E_{D_{\alpha}}^R\left(v^n\right) = E_{\theta_{\alpha}\left(D_{\alpha}\right)}^R\left(v^n\right) = 0$ for every $n \neq 0$. Thus, thanks to the trace property \eqref{Eq10}, Corollary \ref{Cor5-3} shows that $M_{\alpha}'\cap M_{\alpha}^{\omega} \subseteq D_{\alpha}^{\omega} \cap C_{\alpha}^{\omega}$. Here, we need two simple lemmas (which will be used not only here but also later too), and the former has been probably known in the context of orthogonal pairs due to Popa \cite{popa:jot83}. The proofs are both straightforward so that the details are left to the reader. 

\begin{lemma}\label{Lem6-1} Let $N$ be a von Neumann algebra, and let $A$ and $B$ be its von Neumann subalgebras. Suppose that there are faithful normal conditional expectations $E_A : N \rightarrow A$, $E_B : N \rightarrow  B$ and $E_{A\cap B} : N \rightarrow A\cap B$ satisfying the condition{\rm :} 
\begin{equation}\label{Eq11}
E_A\circ E_B = E_{A\cap B}. 
\end{equation}
Then, we have $A^{\omega} \cap B^{\omega} = \left(A\cap B\right)^{\omega}$ inside $N^{\omega}$. 
\end{lemma} 

\begin{lemma}\label{Lem6-2} Let $\gamma : G \rightarrow \mathrm{Aut}(P)$ be an action of a discrete group $G$ on a $\sigma$-finite von Neumann algebra $P$, and assume that $\phi$ is a faithful normal invariant state on $P$ under the action $\gamma$. Let $N :=P\rtimes_{\gamma}G$ be the crossed-product with the canonical unitary representation $\lambda^{\gamma} : G \rightarrow N$ and the canonical conditional expectation $E_P^N : N\rightarrow P$. Let $A$ be a von Neumann subalgebra of $P$, which is the range of a faithful normal conditional expectation $E_A^P : P \rightarrow A$. Let us also choose a von Neumann subalgebra $B$ of $\lambda^{\gamma}(G)''$, and denote by $E_B^{\lambda^{\gamma}(G)''} : \lambda^{\gamma}(G)'' \rightarrow B$ the conditional expectation that preserves the canonical tracial state $\tau_G$ on $\lambda^{\gamma}(G)''$. Then, we have the following{\rm :}
\begin{itemize} 
\item[(i)] The Fubini map $\phi\otimes\mathrm{Id}_{B\left(\ell^2(G)\right)}$ gives a conditional expectation $E_{\lambda^{\gamma}(G)''}^N : N \rightarrow \lambda^{\gamma}(G)''$ as the restriction to $N$\rm; 
\item[(ii)] The faithful normal conditional expectations 
\begin{equation*}
E_A :=E_A^P \circ E_P^N : N \rightarrow A, \quad   
E_B :=E_B^{\lambda^{\gamma}(G)''}\circ E_{\lambda^{\gamma}(G)''}^N : N \rightarrow B
\end{equation*}
satisfy the condition \eqref{Eq11} with $A \cap B = {\mathbf C}1$ so that $A^{\omega} \cap B^{\omega} = {\mathbf C}1$ inside $N^{\omega}$ due to Lemma \ref{Lem6-1}.
\end{itemize}
\end{lemma} 

It is known that $R = C\left({\mathbb T}_{\alpha}\right)''$ can be identified with the crossed-product $L^{\infty}\left({\mathbb T}\right)\rtimes_{\gamma_{\alpha}}{\mathbb Z}$ by the correspondence $u_{\alpha} \leftrightarrow \pi_{\gamma_{\alpha}}\left(\mathrm{id}\right)$, $v_{\alpha} \leftrightarrow \lambda^{\gamma_{\alpha}}(1)$, where $\gamma_{\alpha}$ is the action induced from the group rotation $\zeta \in {\mathbb T} \mapsto e^{2\pi\sqrt{-1}\alpha}\zeta \in {\mathbb T}$ and ``$\mathrm{id}$" denotes the function $\mathrm{id}(\zeta) :=\zeta$, $\zeta \in {\mathbb T}$. Thus, Lemma \ref{Lem6-2} implies that $D_{\alpha}^{\omega} \cap C_{\alpha}^{\omega} = {\mathbf C}1$, and hence $M_{\alpha}'\cap M_{\alpha}^{\omega} = \mathbf{C}1$. Summarizing the discussions so far we conclude 
 
\begin{theorem} The HNN extension $M_{\alpha}$ obtained from 
the non-commutative torus $C({\mathbb T}^2_{\alpha})$ with irrational $\alpha \in [0,1)$ in the above manner always becomes a full type {\rm II}$_1$ factor generated by two Haar unitaries. 
\end{theorem}  

\subsubsection{Remark and Question.} A similarity (in some sense) between free entropy dimensions (see \cite{voiculescu:bulllondon2002}) and costs of equivalence relations (see \cite{gaboriau:invent2000}) gives us the question whether or not the (modified) free entropy dimension of $u_{\alpha}$ (or $v_{\alpha}$) and $u(\alpha)$ is $1$. Moreover, we do not know  whether or not $M_{\alpha}$ depends on the choice of $\alpha$.  

\subsection{HNN Extensions Associated with Tensor Product Algebras.} 
\subsubsection{General Setting.} For given   
\begin{itemize} 
\item $\sigma$-finite von Neumann algebras $N_1$, $N_2$, $N_3$; 
\item two isomorphic von Neumann subalgebras 
$D_1$ ($\subseteq N_1$), $D_2$ ($\subseteq N_2$) 
with a surjective (i.e., automatically, normal) $*$-isomorphism $\theta_{21} : D_1 \rightarrow D_2$;  
\item a von Neumann subalgebra $D_3$ ($\subseteq N_3$) 
with an automorphism $\theta_3 \in \mathrm{Aut}\left(D_3\right)$; 
\item faithful normal conditional expectations 
$E_{D_1}^{N_1} : N_1 \rightarrow D_1$, 
$E_{D_2}^{N_2} : N_2 \rightarrow D_2$, 
$E_{D_3}^{N_3} : N_3 \rightarrow D_3$,
\end{itemize} 
we set 
\begin{equation*}
N :=N_1\otimes N_2\otimes N_3, \quad  
D :=D_1\otimes{\mathbf C}1\otimes D_3,
\end{equation*}
and define the normal $*$-isomorphism $\theta : D \rightarrow N$ by 
\begin{equation*}
\theta : d_1\otimes 1\otimes d_3 \in D \mapsto 
1\otimes\theta_{21}\left(d_1\right)\otimes\theta_3\left(d_3\right)
\in {\mathbf C}1\otimes  D_2\otimes D_3 \subseteq N
\end{equation*}
(i.e., the interchange between the 1st and 2nd tensor components by the surjective $*$-isomorphism $\theta_{21} : D_1 \rightarrow D_2$ together with the $*$-automorphism $\theta_3$ on the 3rd tensor component $D_3$). Here, we mention that the third component $N_3$ allows a type III$_0$ example, see \S\S\S6.2.6, Case 2. For given faithful normal states $\varphi_1$, $\varphi_2$ on 
$D_1$, $D_2$, respectively, we define  
\begin{equation*}
E_D^N :=
E_{D_1}^{N_1}\otimes\left(\varphi_2\circ E_{D_2}^{N_2}\right)\otimes E_{D_3}^{N_3}, \quad  
E_{\theta(D)}^N :=
\left(\varphi_1\circ E_{D_1}^{N_1}\right)\otimes E_{D_2}^{N_2}\otimes E_{D_3}^{N_3}. 
\end{equation*}
Let us construct and investigate the HNN extension: 
\begin{equation*}
\left(M, E_N^M\right) = 
\left(N, E_D^N\right) \underset{D}{\bigstar} \left(\theta, E_{\theta(D)}^N\right)
\end{equation*}
with stable unitary $u(\theta)$, where the embedding map of $N$ into $M$ is omitted as before. 

\subsubsection{Assumption.} In what follows, let us assume the condition:
There are two unitaries 
$v_1 \in \left(N_1\right)_{\varphi_1\circ E_{D_1}^{N_1}}$, 
$v_2 \in \left(N_2\right)_{\varphi_2\circ E_{D_2}^{N_2}}$ satisfying
\begin{equation*}
\varphi_1\circ E_{D_1}^{N_1}\left(v_1^n\right) = 
\varphi_2\circ E_{D_2}^{N_2}\left(v_2^n\right) = 0 
\end{equation*}
as long as $n \neq 0$. 

\subsubsection{Consequences from Corollary \ref{Cor5-3}.} 
With letting $v := v_1\otimes v_2 \otimes 1$ it is plain to see that $E_D^N\left(v^n\right) = E_{\theta(D)}^N\left(v^n\right) = 0$ as long as $n\neq 0$. For a fixed faithful normal state $\varphi_3$ on $D_3$, we define the two states $\varphi$, $\varphi_{\theta}$ on $D$ in such a way that    
\begin{equation*}
\varphi\left(d_1\otimes d_3\right) :=
\varphi_1\otimes\varphi_3\left(d_1\otimes d_3\right), \quad 
\varphi_{\theta}\left(d_1\otimes d_3\right) :=
\varphi_2\otimes\varphi_3\left(\theta_{21}\left(d_1\right)\otimes\theta_3\left(d_3\right)\right) 
\end{equation*}
for each $d_1\otimes d_3 \in D_1\otimes D_3 \cong D_1\otimes{\mathbf C}1\otimes D_3 = D$ (with the identification $d_1\otimes d_3 \leftrightarrow d_1\otimes 1 \otimes d_3$), and hence  
{\allowdisplaybreaks 
\begin{equation*} 
\varphi\circ E_D^N = \left(\varphi_1\circ E_{D_1}^{N_1}\right)\otimes\left(\varphi_2\circ E_{D_2}^{N_2}\right)\otimes\left(\varphi_3\circ E_{D_3}^{N_3}\right) = \varphi_{\theta}\circ\theta^{-1}\circ E_{\theta(D)}^N. 
\end{equation*}
}Thus the unitary $v$ belongs to both the centralizers $N_{\varphi\circ E_D^N} = N_{\varphi_{\theta}\circ\theta^{-1}\circ E_{\theta(D)}^N}$. Therefore, Corollary \ref{Cor5-3} implies that  
\begin{eqnarray} 
{\mathcal Z}\left(M\right) &=&  
\left\{ x \in D \cap \theta(D) \cap N' : \theta(x) = x \right\},  \label{Eq12}\\
{\mathcal Z}\left(\widetilde{M}\right) &=&  
\left\{ x \in \widetilde{D} \cap \widetilde{\theta}\left(\widetilde{D}\right) \cap \widetilde{N}' : 
\widetilde{\theta}(x) = x \right\}, \label{Eq13}\\
M'\cap M^{\omega} &=&  
\left\{ x \in D^{\omega} \cap \theta(D)^{\omega} \cap N' : 
\theta^{\omega}(x) = x \right\}. \label{Eq14}
\end{eqnarray}

\subsubsection{Factoriality and Flow of Weights.} To investigate the factoriality and the flow of weights of $M$, it suffices to determine the right-hand sides of \eqref{Eq12} and \eqref{Eq13} explicitly. Thanks to \cite[Vol.I; Theorem 5.9, Corollary 5.10 in Chap.~IV]{takesaki:book}, we see that 
{\allowdisplaybreaks 
\begin{align*}
D\cap\theta(D)\cap N' 
&= {\mathbf C}1\otimes{\mathbf C}1\otimes\left(D_3\cap N_3'\right), 
\end{align*}
}and hence, by the definition of $\theta$ we get 
\begin{equation*}
\left\{ x \in D \cap \theta(D) \cap N'\ :\ \theta(x) = x \right\} \cong D_3^{\theta_3}\cap N_3'. 
\end{equation*}
With letting $\displaystyle{\sigma_t :=\sigma_t^{\varphi\circ E_D^N} =  
\sigma_t^{\varphi_1\circ E_{D_1}^{N_1}}\otimes
\sigma_t^{\varphi_2\circ E_{D_2}^{N_2}}\otimes
\sigma_t^{\varphi_3\circ E_{D_3}^{N_3}}}$, 
the continuous cores $\widetilde{N} \supseteq \widetilde{D},\ \widetilde{\theta}\left(\widetilde{D}\right) = \widetilde{\theta(D)}$ are captured as the simultaneous fixed-point algebras:  
{\allowdisplaybreaks 
\begin{align*}
\widetilde{N} &= \left(N_1\otimes N_2\otimes N_3\right)\rtimes_{\sigma}{\mathbf R} = 
\left(\left(N_1\otimes N_2\otimes N_3\right)\otimes B\left(L^2\left({\mathbf R}\right)\right)\right)
^{\sigma_t\otimes\mathrm{Ad}\lambda_{-t}}, \\
\widetilde{D} &= \left(D_1\otimes{\mathbf C}1\otimes D_3\right)\rtimes_{\sigma}{\mathbf R} = 
\left(\left(D_1\otimes{\mathbf C}1\otimes D_3\right)\otimes 
B\left(L^2\left({\mathbf R}\right)\right)\right)
^{\sigma_t\otimes\mathrm{Ad}\lambda_{-t}}, \\
\widetilde{\theta}\left(\widetilde{D}\right) &= 
\left({\mathbf C}1\otimes D_2\otimes D_3\right)\rtimes_{\sigma}{\mathbf R} = 
 \left(\left({\mathbf C}1\otimes D_2\otimes D_3\right)\otimes 
 B\left(L^2\left({\mathbf R}\right)\right)\right)
^{\sigma_t\otimes\mathrm{Ad}\lambda_{-t}}
\end{align*}
}by the Takesaki duality theorem (see \cite[Vol.II; Theorem 2.3 in Ch.~X]{takesaki:book}). Thus, we have
{\allowdisplaybreaks 
\begin{align*} 
\widetilde{D}\cap\widetilde{\theta}\left(\widetilde{D}\right) 
&= 
{\mathbf C}1\otimes{\mathbf C}1\otimes 
\left(D_3\otimes  B\left(L^2\left({\mathbf R}\right)\right)
\right)^{\sigma_t^{\varphi_3\circ E_{D_3}^{N_3}}\otimes\mathrm{Ad}\lambda_{-t}} \\ 
&=  
{\mathbf C}1\otimes{\mathbf C}1\otimes 
\left(D_3\rtimes_{\sigma^{\varphi_3\circ E_{D_3}^{N_3}}}{\mathbf R}\right),  
\end{align*}
}and then 
\begin{equation*} 
\left\{ x \in \widetilde{D} \cap \widetilde{\theta}\left(\widetilde{D}\right) \cap \widetilde{N}'\ :\ 
\widetilde{\theta}(x) = x \right\}  
= 
\left({\mathbf C}1\otimes{\mathbf C}1\otimes 
\left(D_3\rtimes_{\sigma^{\varphi_3\circ E_{D_3}^{N_3}}}
{\mathbf R}\right)^{\widetilde{\theta}_3}\right)
\cap\widetilde{N}'  
\end{equation*}
with the canonical extension $\widetilde{\theta_3}$ of $\theta_3$. Summing up what we have done, we conclude 

\begin{theorem}\label{Thm6-4} Under Assumption 6.2.2, we have 
\begin{equation*}
{\mathcal Z}\left(M\right) \cong D_3^{\theta_3}\cap N_3' = D_3^{\theta_3}\cap{\mathcal Z}\left(N_3\right), 
\end{equation*}
and in particular, if $N_3$ is a factor or $D_3^{\theta_3}$ is 
the trivial algebra ${\mathbf C}1$, then so is $M$. Moreover, 
we also have 
\begin{equation*}
{\mathcal Z}\left(\widetilde{M}\right) 
=  \left({\mathbf C}1\otimes{\mathbf C}1\otimes\left(\widetilde{D_3}^{\widetilde{\theta_3}}\cap{\mathcal Z}\left(\widetilde{N_3}\right)\right)\right) 
\cap {\mathcal Z}\left(\widetilde{N}\right),   
\end{equation*}
where $\widetilde{N_3} \supseteq \widetilde{D_3}$ is the inclusion of the continuous cores of $N_3 \supseteq D_3$ determined by $E_{D_3}^{N_3}$ and $\widetilde{\theta_3}$ is the canonical extension of $\theta_3$.
\end{theorem} 

\subsubsection{Fullness.} To determine whether or not $M$ is full, we first have to determine the right-hand side of \eqref{Eq14}.  To do so, we will use Lemma \ref{Lem6-1}. For every $x \in N$, we have $E_D^N\left(E_{\theta(D)}^N\left(x\right)\right) = \left(\left(\varphi_1\circ E_{D_1}^{N_1}\right)\otimes\left(\varphi_2\circ E_{D_2}^{N_2}\right)\otimes E_{D_3}^{N_3}\right)(x) 
\in {\mathbf C}1\otimes{\mathbf C}1\otimes D_3$.
Note that $D\cap \theta(D) = {\mathbf C}1\otimes{\mathbf C}1\otimes D_3$ thanks to \cite[Vol.II; Corollary 5.10 in Chap.~IV]{takesaki:book}, and hence the condition in Lemma \ref{Lem6-1} holds. Therefore we get: 

\begin{theorem}\label{Thm6-5} Under Assumption 6.2.2, we have 
\begin{equation*}
M_{\omega} \subseteq 
M'\cap M^{\omega} \cong \left(D_3^{\omega}\right)^{\theta_3^{\omega}}\cap N_3'.
\end{equation*}
In particular, if $N_3$ is the trivial algebra, then $M$ always becomes a full factor. Also, if the right-hand side sit in the asymptotic centralizer $\left(N_3\right)_{\omega}$, then the first inclusion relation would become the identity.   
\end{theorem} 

\subsubsection{More Concrete Cases.} We should first remark that there is a variety of concrete examples which satisfy Assumption 6.2.2. In fact, if a given pair $(L, \psi)$ of von Neumann algebra and faithful normal state had the non-atomic centralizer $L_{\psi}$, one would be able to find a Haar unitary in $L_{\psi}$ with respect to $\psi$. Thus, in what follows, we may assume that the given two quartets $\left(N_1 \supseteq D_1, E_{D_1}^{N_1}, \varphi_1\right)$,  $\left(N_2 \supseteq D_2, E_{D_2}^{N_2}, \varphi_2\right)$ satisfy Assumption 6.2.2.  

\medskip\noindent
{\bf Case 1.} Assume that $N_3$ is the trivial algebra, that is, no presence of the triple $N_3 \supseteq D_3, \theta_3 \in \mathrm{Aut}\left(D_3\right)$ in our initial data $N \supseteq D$, $\theta : D \rightarrow N$. Then, Theorem \ref{Thm6-4} and Theorem \ref{Thm6-5} say that $M$ is a full factor, and the center of its continuous core is computed as follows. 
\begin{equation}\label{Eq15}
{\mathcal Z}\left(\widetilde{M}\right) = \left({\mathbf C}\rtimes_{\sigma}{\mathbf R}\right)\cap{\mathcal Z}\left(\left(N_1\otimes N_2\right)\rtimes_{\sigma}{\mathbf R}\right) 
\end{equation} 
with $\sigma_t = \sigma_t^{\varphi_1\circ E_{D_1}^{N_1}}\otimes\sigma_t^{\varphi_2\circ E_{D_2}^{N_2}}$. The above \eqref{Eq15}, in particular, shows that the flow of weights of $M$ is a factor flow of the translation on the real line ${\mathbf R}$ so that $M$ does never become of type III$_0$. To find the exact number ``$\lambda$" in the III$_{\lambda}$-classification, one must determine the right-hand side of \eqref{Eq15} in more concrete form or the T-set $T(M)$ very explicitly, both of which seem somewhat delicate tasks except for several simple cases. We will next illustrate how the T-set can be determined in one of such simple cases. Assume further that two states $\varphi_1$, $\varphi_2$ are chosen so that  $\varphi_2\circ\theta_{21} = \varphi_1$. Then, Theorem \ref{Thm4-1} implies that $\sigma_t^{\varphi\circ E_D^N\circ E_N^M}\left(u(\theta)\right) = u(\theta)$, and hence by Theorem \ref{Thm5-2} we have $\left(M_{\varphi\circ E_D^N\circ E_N^M}\right)' \cap M \subseteq \left\langle u(\theta) \right\rangle' \cap N$. Here, we do the same argument as in the first part of the proof of Corollary \ref{Cor5-3} and get $\left(M_{\varphi\circ E_D^N\circ E_N^M}\right)' \cap M = {\mathbf C}1$ since $D \cap \theta(D) = {\mathbf C}1$. This computation implies that 
\begin{equation*} 
T(M) = \left\{ t \in {\mathbf R}\ :\ \sigma_t^{\varphi_1\circ E_{D_1}^{N_1}} = \mathrm{Id} = \sigma_t^{\varphi_2\circ E_{D_2}^{N_2}} \right\}, 
\end{equation*}
whose right-hand side is computable when all the initial data are given concretely.     

\medskip\noindent
{\bf Case 2.} Assume that $D_3$ is a non-trivial algebra, while $\theta_3 \in \mathrm{Aut}(D_3)$ is assumed to be ergodic or $N_3$ to be a factor. Then, Theorem \ref{Thm6-4} says that $M$ is a factor. We further assume that both $\varphi_1\circ E_{D_1}^{N_1}$ and $\varphi_2\circ E_{D_2}^{N_2}$ are traces. In this case, it is plain to see, by Theorem \ref{Thm6-4}, that 
\begin{equation*}
{\mathcal Z}\left(\widetilde{M}\right) \cong \widetilde{D_3}^{\widetilde{\theta_3}}\cap{\mathcal Z}\left(\widetilde{N_3}\right)
\end{equation*} 
so that the type of $M$ is completely determined from the data $N_3 \supseteq D_3$, $\theta_3 \in \mathrm{Aut}\left(D_3\right)$, and $M$ can be of type III$_0$ in this case. Instead of assuming that $\varphi_1\circ E_{D_1}^{N_1}$ and $\varphi_2\circ E_{D_2}^{N_2}$ are traces, we will next impose the extra assumption that the states $\varphi_1$, $\varphi_2$ and $\varphi_3$ are chosen so that $\varphi_2\circ\theta_{21} = \varphi_1$ and $\varphi_3\circ\theta_3 = \varphi_3$. Then, Theorem \ref{Thm4-1} implies that $\sigma_t^{\varphi\circ E_D^N\circ E_N^M}\left(u(\theta)\right) = u(\theta)$, and hence by the same argument as in the final part of Case 1, we get $\left(M_{\varphi\circ E_D^N\circ E_N^M}\right)'\cap M = {\mathbf C}1$. Thus, $M$ is not of type III$_0$, and the T-set $T(M)$ is computed as 
\begin{equation*}
T(M) = \left\{ t \in {\mathbf R}\ :\ \sigma_t^{\varphi_1\circ E_{D_1}^{N_1}} = \sigma_t^{\varphi_2\circ E_{D_2}^{N_2}} = \sigma_t^{\varphi_3\circ E_{D_3}^{N_3}}= \mathrm{Id}\right\}.  
\end{equation*}
In all the cases treated in Case 2, it seems difficult to determine the asymptotic centralizer $M_{\omega}$ (or whether $M$ is full or not) except for the case that $D_3$ is a finite von Neumann algebra because it is non-trivial in general whether $M_{\omega} = M'\cap M^{\omega}$ or not.   
 
\subsection{HNN extensions arising from pairs of regular and singular MASAs}  
\subsubsection{General Setting.} 
Let $Q$ be a $\sigma$-finite von Neumann algebra with a faithful normal state $\varphi_Q$, and $G$ be a countably infinite discrete group. We then construct the infinite tensor product over $G$: 
\begin{equation*}
\left(P, \varphi_P\right) :=
\bigotimes_{g \in G} \left(Q, \varphi_Q\right)_g, 
\end{equation*}
where the $\left(Q,\varphi_Q\right)_g$'s are copies of $\left(Q,\varphi_Q\right)$. The canonical embedding map of $Q$ onto the $g$th tensor component in $P$ is denoted by $\iota_g$. By the construction, for distinct elements $g_1,\dots,g_n \in G$ and for $x_1,\dots,x_n \in Q$, the operators $\iota_{g_1}\left(x_1\right),\dots,\iota_{g_n}\left(x_n\right)$ mutually commute and $\varphi_P\left(\iota_{g_1}\left(x_1\right)\cdots\iota_{g_n}\left(x_n\right)\right) = \varphi_Q\left(x_1\right)\cdots\varphi_Q\left(x_n\right)$. Let us denote by $\gamma$ the Bernoulli shift action of $G$ on $P$ defined in such a way that $\gamma_h\left(\iota_g\left(x\right)\right) :=\iota_{hg}\left(x\right)$ ($x \in Q$, $h, g \in G$), and we will consider the crossed-product $N := P\rtimes_{\gamma}G$ with the canonical unitary representation $\lambda^{\gamma} : G \rightarrow N$ and the canonical conditional expectation $E_P^N : N \rightarrow P$.
Since $\gamma$ is invariant under the product state $\varphi_P$, the Fubini map $\varphi_P\otimes\mathrm{Id}_{B\left(\ell^2(G)\right)}$ gives a faithful normal conditional expectation from $N$ onto $\lambda^{\gamma}(G)''$, see Lemma \ref{Lem6-2}, and it is denoted by $E_{\lambda^{\gamma}(G)''}^N$. Let us choose a subgroup $H$ of $G$ and a von Neumann subalgebra $D$ of $P$ with the $\varphi_P$-preserving conditional expectation $E_D^P : P \rightarrow D$ in such a way that $\left(D, \varphi_P|_D\right) \cong \left(\lambda^{\gamma}(H)'',\tau_G|_{\lambda^{\gamma}(H)''}\right)$ in the state-preserving way. Letting $\tau := \varphi_P|_D$ and $\tau_H := \tau_G|_{\lambda^{\gamma}(H)''}$, we have a surjective $*$-isomorphism $\theta : D \rightarrow \lambda^{\gamma}(H)''$ with the property:  
\begin{equation}\label{Eq16}
\tau = \tau_H\circ\theta.
\end{equation}
Let $E_{\theta(D)}^{\lambda^{\gamma}(G)''} : \lambda^{\gamma}(G)'' \rightarrow \theta(D) = \lambda^{\gamma}(H)''$ be the $\tau_G$-preserving conditional expectation, and set 
\begin{equation*}
E_D^N :=E_D^P\circ E_P^N, \quad E_{\theta(D)}^N :=E_{\theta(D)}^{\lambda^{\gamma}(G)''}\circ E_{\lambda^{\gamma}(G)''}^N,  
\end{equation*}
conditional expectations from $N$ onto $D$ and $\theta(D) = \lambda^{\gamma}(H)''$, respectively. Then we construct the HNN extension 
\begin{equation*}
\left(M,E_N^M\right) = \left(N,E_D^N\right) \underset{D}{\bigstar} \left(\theta, E_{\theta(D)}^N\right)
\end{equation*}
with stable unitary $u(\theta)$, where the embedding map of $N$ into $M$ is omitted as before. 

\subsubsection{Assumptions.} In what follows, we will assume that (i) there is a unitary $u$ in the centralizer $Q_{\varphi_Q}$ such that $\varphi_Q(u) = 0$; (ii) $G$ has an element $g_0 \in G$ of infinite order. 

\subsubsection{Consequences from Corollary \ref{Cor5-3}.} Letting $v :=\iota_{e}(u)\lambda^{\gamma}\left(g_0\right) \in N$, we have, for every $n \in {\mathbb N}$,   
{\allowdisplaybreaks 
\begin{align*}  
E_D^N\left(v^n\right) &= E_D^P\left(\left(\iota_{e}\left(u\right)\cdots\iota_{g_0^{n-1}}\left(u\right)\right) E_P^N\left(\lambda^{\gamma}\left(g_0^n\right)\right)\right) = 0; \\
E_{\theta(D)}^N\left(v^n\right) 
&= \varphi_P\left(\iota_{e}\left(u\right)\cdots\iota_{g_0^{n-1}}\left(u\right)\right) E_{\theta(D)}^{L(G)}\left(\lambda^{\gamma}\left(g_0^n\right)\right) 
= 0. 
\end{align*}
Note that  
\begin{equation}\label{Eq17}
\tau\circ E_D^N = \varphi_P\circ E_P^N = \tau_G\circ E_{L(G)}^N = \tau_H\circ E_{\theta(D)}^N, 
\end{equation} 
and hence by \eqref{Eq16}
\begin{equation}\label{Eq18}
\tau\circ E_D^N = \tau_H\circ E_{\theta(D)}^N = 
\tau\circ\theta^{-1}\circ E_{\theta(D)}^N.
\end{equation}
By the above computation \eqref{Eq17} we get $\sigma_t^{\tau\circ E_D^N}(v) =v$ so that the unitary $v$ is in the centralizer $N_{\tau\circ E_D^N} = N_{\tau\circ\theta^{-1}\circ E_{\theta(D)}^N}$. Thus, Corollary \ref{Cor5-3} shows that $M_{\omega} \subseteq M' \cap M^{\omega} \subseteq D^{\omega}\cap\theta^{\omega}\left(D^{\omega}\right)$. Furthermore, Lemma \ref{Lem6-1} and Lemma \ref{Lem6-2} show that $D^{\omega}\cap\theta^{\omega}\left(D^{\omega}\right) = \mathbf{C}1$. Therefore, we conclude 

\begin{theorem}\label{Thm6-6}
Under Assumptions 6.3.2, the HNN extension $M$ is a full factor. 
\end{theorem}  

\subsubsection{Modular Automorphisms and Type Classification.} Thanks to the above \eqref{Eq18} together with Theorem \ref{Thm4-1}, we observe  
\begin{equation}\label{Eq19}
\sigma_t^{\tau\circ E_D^N\circ E_N^M}\left(u(\theta)\right) = u(\theta) \left[D\tau\circ\theta^{-1}\circ E_{\theta(D)}^N : D\tau\circ E_D^N\right]_t = u(\theta),  
\end{equation} 
and hence the type classification of $M$ is the same as that of the crossed-product $N$. In fact, Theorem \ref{Thm5-2} together with \eqref{Eq19} implies that 
\begin{equation*}
\left(M_{\tau\circ E_D^N\circ E_N^M}\right)' \cap M \subseteq 
\left\langle v, u(\theta) \right\rangle' \cap M \subseteq N
\end{equation*}
so that $\left(M_{\tau\circ E_D^N\circ E_N^M}\right)' \cap M  \subseteq \left\langle u(\theta) \right\rangle' \cap N$. Now, by the same argument as in the first part of the proof of Corollary \ref{Cor5-3}, we see that the right-hand side of the above sits in $D \cap \theta(D)$ being the trivial algebra ${\mathbf C}1$ by Lemma \ref{Lem6-2}. Therefore, we have the trivial relative commutant property $\left(M_{\tau\circ E_D^N\circ E_N^M}\right)' \cap M = {\mathbf C}1$. This says that the T-set $T(M)$ is enough to determine the type of $M$ and in particular that $M$ can never become of type II$_{\infty}$ nor type III$_0$ (see e.g.~the discussions given in \cite[p.377-388]{ueda:pacific}). Moreover, the relative commutant property implies that 
\begin{equation*}
T(M) = \left\{ t \in {\mathbf R}\ :\ \sigma_t^{\tau\circ E_D^N\circ E_N^M} = \mathrm{Id} \right\} = \left\{ t \in {\mathbf R}\ :\ \sigma_t^{\tau\circ E_D^N} = \mathrm{Id} \right\}.  
\end{equation*}
Here, the second equality simply comes from \eqref{Eq19}. Since the product state $\varphi_P$ is invariant under the Bernoulli shift action $\gamma$, we observe that  
\begin{equation*}
N_{\varphi_P\circ E_P^N} = \left\langle P_{\varphi_P},\ \lambda^{\gamma}(G) \right\rangle'' =  
P_{\varphi_P}\rtimes_{\gamma}G
\end{equation*}
(see e.g.~the discussion on the free analog of Connes-St{\o}rmer's Bernoulli shifts in \cite{hiaiueda:jot}), and this algebra is clearly a type II$_1$ factor. Hence, by \eqref{Eq17}  
\begin{equation*}
\sigma_t^{\tau\circ E_D^N} = \mathrm{Id} \Longleftrightarrow \sigma_t^{\varphi_P} = \mathrm{Id} \Longleftrightarrow \sigma_t^{\varphi_Q} = \mathrm{Id}. 
\end{equation*}
Therefore, we conclude 

\begin{theorem}\label{Thm6-7} Under Assumptions 6.3.2, the HNN extension $M$ is of type II$_1$ or type III$_{\lambda}$ with $\lambda \neq 0$, and the type is completely determined by the T-set. The T-set is computed as follows. 
\begin{equation}\label{Eq20}
T(M) = T(N) = \left\{ t \in {\mathbf R}\ :\ \sigma_t^{\varphi_Q} = \mathrm{Id} \right\}.  
\end{equation}
\end{theorem}
 
\subsubsection{Concrete Cases{\rm :} Regular vs Singular MASAs} We will give type II$_1$ and type III$_{\lambda}$ ($0 < \lambda \leq 1$) concrete examples of HNN extensions $M = N \underset{D}{\bigstar} \theta$ such that $D$ is a regular MASA in $N$ while $\theta(D)$ a singular MASA in $N$ (so that $\theta$ cannot be extended to any $*$-automorphism of $N$). In what follows, all von Neumann algebras that we will deal with are assumed to have separable preduals. Let $D_0$ be a von Neumann subalgebra of the centralizer $Q_{\varphi_Q}$, and denote by $D$ the von Neumann subalgebra generated by the $\iota_g\left(D_0\right)$'s, that is, $(D,\tau)$ is the infinite tensor product of 
$\left(D_0,\varphi_Q|_{D_0}\right)$ over $G$. Since $D$ sits in the centralizer $P_{\varphi_P}$, the $\varphi_P$-preserving conditional expectation $E_D^P : P \rightarrow D$ exists.   

The next lemma seems a folklore, and the main part of its proof is actually the same as that of showing that non-commutative Bernoulli shifts are free actions. The details are left to the reader. 

\begin{lemma}\label{Lem6-8} 
If $D_0$ is a MASA in $Q$, then so is $D$ in $N$. Furthermore, if $D_0$ is regular in $Q$, then so is $D$ in $N$.  
\end{lemma}

Let us give the pair $\left(Q, \varphi_Q\right)$ more concretely. 

\medskip\noindent
{\bf Case 1.} We discuss the case of $Q = D_0$, and treat the following two cases in common: (1) $Q= D_0$ is the abelian von Neumann algebra of finite dimension greater than $2$; (2) $Q = D_0$ is a diffuse abelian von Neumann algebra. In the case (2), we will be able to treat an arbitrary faithful normal state because one can easily find a Haar unitary with respect to the given state $\varphi_Q$, while in the case (1), the given state $\varphi_Q$ should be constructed from the equal probability vector for the requirement of Assumptions 6.3.2. Assume further that $G = H ={\mathbb Z}$ and $g_0 = 1 \in {\mathbb Z}$. Then, $\left(D=P,\varphi_P\right) \cong \left(L^{\infty}[0,1], \text{Lebesgue measure}\right) \cong \left(\lambda^{\gamma}\left({\mathbb Z}\right)'', \tau_{\mathbb Z}\right)$, and hence a surjective normal $*$-isomorphism $\theta : D \rightarrow \lambda^{\gamma}\left({\mathbb Z}\right)''$ with the property \eqref{Eq16} exists. Thus, the given data $Q = D_0 = {\mathbf C}^n$, $\varphi_Q$, $u$, $G=H={\mathbb Z}$, $g_0=1$ can be treated in the framework that we have worked out in this subsection. Hence,  by Theorem \ref{Thm6-6} and Theorem \ref{Thm6-7}, $M$ is a full factor of type II$_1$. The von Neumann subalgebra $D$ is clearly (from the setting here) a regular MASA in the base algebra $N$, while it is known that $\theta(D) = \lambda^{\gamma}\left({\mathbb Z}\right)''$ is a singular MASA in $N$ thanks to \cite{nielsen:jfa70} (also see \cite[Theorem 2.1]{neshveyevstormer:jfa02} for a recent elegant proof). It is clear that $M$ is generated by $D$ and the stable unitary, and also that $D$ has a Haar unitary generator. Hence, $M$ is generated by two Haar unitaries. 
  
\medskip\noindent
{\bf Case 2.} Let us assume that $Q$ is a non-type I factor with separable predual and $D_0$ is a Cartan subalgebra (i.e., in particular, a regular MASA) in $Q$ with the unique conditional expectation $E_{D_0}^Q : Q \rightarrow D_0$. Thanks to \cite[Lemma 4.2]{ueda:pacific}, one can find a faithful normal state $\psi$ on $D_0$ and a unitary $u \in Q_{\psi\circ E_{D_0}^Q}$ satisfying that  $E_{D_0}^Q\left(u^n\right) = 0$ for every $n \neq 0$. With letting $\varphi_Q :=\psi\circ E_{D_0}^Q$, the triple $\left(Q,\varphi_Q,u\right)$ satisfies the condition (i) in Assumptions 6.3.2. Note here that the subalgebra $D_0$ sits in the centralizer $Q_{\varphi_Q}$. As in Case 1, we assume that $G = H = {\mathbb Z}$ and $g_0 = 1 \in {\mathbb Z}$, and get a surjective $*$-isomorphism $\theta : D \rightarrow \lambda^{\gamma}\left({\mathbb Z}\right)''$ with the property \eqref{Eq16}. Hence, these given data $Q \supseteq D_0$, $\varphi_Q$, $u$, $G = H = {\mathbb Z}$, $g_0 = 1$ can be treated in our framework. Hence, Theorem \ref{Thm6-6} and Theorem \ref{Thm6-7} say that $M$ is a full factor, not of type III$_0$, and the T-set is computed as \eqref{Eq20}. Note that the T-set $T(M)$ does not coincide, in general, with the T-set $T(Q)$ of the initially given factor $Q$ since the right-hand side of \eqref{Eq20} does depend upon the choice of the state $\psi$ on $D_0$. We would like to emphasize that this example of HNN extension can be regarded as a type III$_{\lambda}$ version of those given in Case 1 when the triple $\left(Q \supseteq D_0, \varphi_Q\right)$ is suitably chosen. In fact, as in Case 1, Lemma \ref{Lem6-8} implies that $D$ becomes a regular MASA not only in $P$ but also even in the base algebra $N = P\rtimes_{\gamma}{\mathbb Z}$, while $\theta(D) = \lambda^{\gamma}\left({\mathbb Z}\right)''$ is a singular MASA in $N$ by \cite[Theorem 2.1 and its remark]{neshveyevstormer:jfa02}.  
Finally, it is easy to give, in this setup, a concrete $\left(Q,\varphi_Q\right)$ in such a way that of the T-set $T(M)$ is computable.     

\section{Reduced HNN Extensions of $C^*$-Algebras}

\subsection{Preliminaries on Reduced Free Products with Amalgamations} Because of the same reason as in the von Neumann algebra case, we need to review reduced amalgamated free products with special emphasis on the r\^{o}le of embedding maps of common amalgamated $C^*$-algebras. 

Let  $C$, $A_s$ ($s \in S$, an index set)  be unital $C^*$-algebras, and we have a unital $*$-isomorphism $\iota_s : C \rightarrow A_s$ for each $s \in S$. Suppose that the $C^*$-subalgebra $\iota_s(C)$ of $A_s$ is the range of a conditional expectation $E_s : A_s \rightarrow \iota_s(C)$ for every $s \in S$. For each $s \in S$, let $X_s$ be the separation and completion of $A_s$ with respect to the pre-norm $a \in A_s \mapsto \Vert E_s\left(a^* a\right)\Vert^{1/2}$ with the canonical map $\eta_s : A_s \rightarrow X_s$. The Banach space $X_s$ is equipped with the $C$-valued inner product $\langle\ \cdot\ |\ \cdot\ \rangle_C$ and the right action of $A_s$ defined in such a way that $\left\langle\eta_s(x)|\eta_s(y)\right\rangle_C := E_s\left(x^* y\right)$ and $\eta_s(x)\cdot a := \eta_s\left(xa\right)$ for each $x, y, a \in A_s$. One can also define the left action of $A_s$, as a $*$-homomorphism into the adjointable operators $B\left(X_s\right)$ on the Hilbert right $A_s$-module $X_s$, defined by $a\cdot\eta_s(x) := \eta_s(ax)$ for each $a, x \in A_s$. This left action of $A_s$ is usually called the GNS representation associated with $E_s$.  We can regard $X_s$ as a $C$-$C$ bimodule (via $\iota_s$) by restricting both the left and right actions of $A_s$ to the subalgebra $\iota_s(C)$. The restriction of $\eta_s$ to $\iota_s(C)$ is clearly injective, and $\eta_s\left(\iota_s(C)\right)$ becomes a complimented closed sub-bimodule of $X_s$ thanks to the decomposition $A_s = \iota_s\left(C\right) + \mathrm{Ker}E_s$. Hence, we have $X_s \cong \iota_s(C)\oplus X_s^{\circ} \cong C\oplus X_s^{\circ}$ as a $C$-$C$ bimodule by the identification $\eta_s\left(\iota_s(c) + a^{\circ}\right) \leftrightarrow \iota_s(c)\oplus\eta_s\left(a^{\circ}\right) \leftrightarrow c\oplus \eta_s\left(a^{\circ}\right)$. 

By the construction in \cite[\S5]{voiculescu:lecturenotes1132} together with the above-mentioned fact on $C^*$-bimodules, we obtain a unital $C^*$-algebra $A$, two kinds of $*$-homomorphisms $\lambda : C \rightarrow A$, $\lambda_s : A_s \rightarrow A$, $s \in S$, and a conditional expectation $E : A \rightarrow \lambda(C)$ satisfying (i) $A$ is generated by the $\lambda_s\left(A_s\right)$'s; (ii) $\lambda_s\circ\iota_s = \lambda$, $s \in S$; (iii) $E\circ\lambda_s = \lambda\circ E_s$, $s \in S$; (iv) the $\lambda_s\left(A_s\right)$'s are free with amalgamation with respect to $E$; (v) if $x \in A$ satisfies $E\left(a^* x^* x a\right) = 0$ for all $a \in A$, then $x = 0$. These five conditions characterize the pair $\left(A, E\right)$ together with $\lambda$ and $\lambda_s$, $s \in S$, completely, see \cite[\S\S5.6]{voiculescu:lecturenotes1132} for details (also see \S2 for more careful explanation on the admissibility of embedding maps of amalgamated algebras). We denote 
\begin{equation*}
\left(A, E\right) = \underset{s\in S}{\bigstar_C} \left(A_s, E_s : \iota_s\right), 
\end{equation*}
and call it the free product of the $A_s$'s with amalgamation over $C$ via $\iota_s$ with respect to the $E_s$'s. 
 
\subsection{Reduced HNN Extensions of $C^*$-Algebras} Let $B$ be a unital $C^*$-algebra and $C$ be a distinguished unital $C^*$-subalgebra with a conditional expectation $E_C^B : B \rightarrow C$. Let us suppose that we have a family $\Theta$ of $*$-isomorphisms $\theta : C \rightarrow  B$ with conditional expectations $E_{\theta\left(C\right)}^B : B \rightarrow \theta\left(C\right)$, $\theta \in \Theta$. In what follows, we will do the same construction as in the von Neumann algebra case. 

Set $\Theta_1 :=\{1 := \mathrm{Id}_D\}\sqcup\Theta$, a disjoint union. 
Define the embedding map 
$\iota_{\Theta} : C\otimes \ell^{\infty}\left(\Theta_1\right) \rightarrow B\otimes B\left(\ell^2\left(\Theta_1\right)\right)$ by 
\begin{equation*}
\iota_{\Theta}\left(x\otimes e_{\theta \theta}\right) :=
\begin{cases} 
x\otimes e_{1 1} & \text{if $\theta = 1$}, \\
\theta(x)\otimes e_{\theta \theta} & \text{if $\theta \in \Theta$}
\end{cases} 
\end{equation*}
for each $x \in C$, where the $e_{\theta_1 \theta_2}$'s denote the canonical matrix unit system in $B\left(\ell^2\left(\Theta_1\right)\right)$, and the conditional expectation $E_{\Theta} : B\otimes B\left(\ell^2\left(\Theta_1\right)\right) \rightarrow \iota_{\Theta}\left(C\otimes\ell^{\infty}\left(\Theta_1\right)\right)$ is defined by 
\begin{equation*}
E_{\Theta} :=\left(\sideset{}{^{\oplus}}\sum_{\theta \in \Theta_1} E_{\theta\left(C\right)}^B\otimes\mathrm{Id}_{{\mathbf C}e_{\theta \theta}}\right)\circ \left(\mathrm{Id}_N\otimes E_{\ell^{\infty}}\right), 
\end{equation*}
where $E_{\ell^{\infty}}$ is the unique conditional expectation from $B\left(\ell^2\left(\Theta_1\right)\right)$ onto $\ell^{\infty}\left(\Theta_1\right)$. Let us denote the inclusion map of $C\otimes\ell^{\infty}\left(\Theta_1\right)$ into $B\otimes B\left(\ell^2\left(\Theta_1\right)\right)$ by $\iota_1$, and define the conditional expectation $E_1 : C\otimes B\left(\ell^2\left(\Theta_1\right)\right) \rightarrow C\otimes\ell^{\infty}\left(\Theta_1\right)$ by 
\begin{equation*}
E_1 :=
\left(E_C^B\otimes\mathrm{Id}_{\ell^{\infty}\left(\Theta_1\right)}\right)\circ\left(\mathrm{Id}_B
\otimes E_{\ell^{\infty}}\right) = E_C^B\otimes E_{\ell^{\infty}}. 
\end{equation*}   
We then construct the reduced free product with amalgamation: 
\begin{equation*}
\left({\mathcal B}, {\mathcal E}\right) = 
\left(B\otimes B\left(\ell^2\left(\Theta_1\right)\right), E_{\Theta} : \iota_{\Theta}\right) 
\underset{C\otimes\ell^{\infty}\left(\Theta_1\right)}{\bigstar} 
\left(B\otimes B\left(\ell^2\left(\Theta_1\right)\right), E_1 : \iota_1\right),     
\end{equation*}
and the embedding maps of $B\otimes B\left(\ell^2\left(\Theta_1\right)\right)$ onto the 1st/2nd free components are denoted by $\lambda_{\Theta}$ and $\lambda_1$, respectively, and the embedding map of $C\otimes\ell^{\infty}\left(\Theta_1\right)$ into ${\mathcal B}$ by $\lambda$ as usual. As in the von Neumann algebra case, we define  
\begin{itemize}
\item $u(\theta) :=\lambda_1\left(e_{1 \theta}\right)\lambda_{\Theta}\left(e_{\theta 1}\right)$ 
with identifying $e_{\theta_1 \theta_2} = 1\otimes e_{\theta_1 \theta_2}$; 
\item the projection $p :=\lambda\left(e_{11}\right) = \lambda_{\Theta}\left(e_{11}\right) \in {\mathcal B}$; 
\item the unital $*$-homomorphism $\pi : B \rightarrow p{\mathcal B}p$  by $\pi(b) :=\lambda_{\Theta}\left(b\otimes e_{11}\right)$ for every $b \in B$. 
\end{itemize}
The partial isometries $u(\theta)$, $\theta \in \Theta$, can be thought of as unitaries in the corner $p{\mathcal B}p$, and by the exactly same way as in the von Neumann algebra case we have the relation $u(\theta)\pi\left(\theta\left(c\right)\right)u(\theta)^* = \pi\left(c\right)$ for each $c \in C, \theta \in \Theta$. Let us denote by $A$ the unital $C^*$-subalgebra of $p{\mathcal B}p$ generated by $\pi\left(B\right)$ and all the $u(\theta)$'s. The restriction of the conditional expectation ${\mathcal E}_{\Theta} : {\mathcal B} \rightarrow \lambda_{\Theta}\left(B\otimes B\left(\ell^2\left(\Theta_1\right)\right)\right)$ conditioned by ${\mathcal E}$ (i.e., ${\mathcal E} = {\mathcal E}\circ{\mathcal E}_{\Theta}$ holds) to $A$ gives rise to a conditional expectation from $A$ onto $\pi(B)$, i.e., $E^A_{\pi(B)} = {\mathcal E}_{\Theta}\big|_A : A \rightarrow \pi(B)$, since $\pi(B) = p\lambda_{\Theta}\left(B\otimes B\left(\ell^2\left(\Theta_1\right)\right)\right)p$. It is easily verified that ${\mathcal E}|_A = \pi\circ  E_C^B\circ\pi^{-1}\circ E_{\pi(B)}^A$. 

\begin{definition} {\rm {\bf (Reduced HNN extensions)}} We call the pair $\left(A, E^A_{\pi(B)}\right)$ constructed so far the reduced HNN extension of $B$ by $\Theta$ with respect to $E_C^B$ and the $E_{\theta(C)}^B$, $\theta \in \Theta$, and  denote it by 
\begin{equation*}
\left(A, E^A_{\pi(B)}\right) = \left(B, E_C^B\right) \underset{C}{\bigstar} 
\left(\Theta, \left\{E_{\theta(C)}^B\right\}_{\theta \in \Theta}
\right). 
\end{equation*}
When no confusion occurs, we will write $A = B\underset{C}{\bigstar}\Theta$ for short. 
\end{definition}  

Not only the notion of reduced words and the conditions (A), (M) are of course valid even in this $C^*$-algebra setting, but also so does the following characterization:  

\begin{proposition} 
The pair $\left(A, E^A_{\pi(B)}\right)$ constructed above satisfies both the conditions {\rm (A)}, {\rm (M)}.  On the other hand, the conditions {\rm (A)}, {\rm (M)} characterize the pair $\left(A, E^A_{\pi(B)}\right)$ completely under the assumptions that {\rm (i)} $\pi(B)$ and the $u(\theta)$'s generate $A$ as $C^*$-algebra and {\rm (ii)} if $x \in A$ satisfies $E^A_{\pi(B)}\left(a^* x^* x a\right) = 0$ for all $a \in A$, then $x =0$. More precisely, the conditions {\rm (A)}, {\rm (M)} determine the conditional expectation $E_{\pi(B)}^A$ completely. 
\end{proposition} 

\begin{remark} {\rm As in the von Neumann algebra case, the following holds: Let $G *_H \theta = \langle G, t : t\theta(h)t^{-1} = h,\ h \in H \rangle$ be an HNN extension of group $G$ with stable letter $t$ by group isomorphism $\theta$ from a subgroup $H$ of $G$ into $G$. The reduced group $C^*$-algebra $C^*_r\left(G *_H \theta\right)$ is identified naturally with the reduced HNN extension of $C^*_r\left(G\right)$ with distinguished unitary $\lambda\left(t\right)$, where all the necessary conditional expectations are chosen as canonical tracial state preserving ones.} 
\end{remark}  

\begin{remark}{\bf [Universal HNN extensions]} {\rm In the $C^*$-algebra setting, there is another choice of HNN extensions, that is, the universal one. The universal {\rm (}or full{\rm )} HNN extension ${\mathfrak A} = B \underset{C}{\bigstar}^{\mathrm{univ}} \Theta$ is defined as the universal $C^*$-algebra generated by $B$ and unitaries $u(\theta)$, $\theta \in \Theta$ with only the relation 
$u(\theta) \theta(c) u(\theta)^* = c$ for every $c \in C,\ \theta \in \Theta$. Thus, there is a $*$-homomorphism from this universal $C^*$-algebra ${\mathfrak A}$ onto $A$ sending $b$ and $u(\theta)$ to $\pi(b)$ and $u(\theta)$, respectively. This means that our construction of reduced HNN extensions can be thought of as a procedure to construct a conditional expectation from ${\mathfrak A}$ onto $B$ by using given $E_C^B$, $E_{\theta(C)}^B$, $\theta \in \Theta$ when $\pi$ is a faithful representation of $B$ {\rm (}this is the case when the given conditional expectations are faithful{\rm )}. The existence of such universal HNN extension can be shown in the same way as in the group case, based on the universal amalgamated free product and the  universal crossed-product constructions.}  
\end{remark} 

\subsection{Embedding of Subsystems in the Framework of Reduced HNN Extensions} Assume that $B_0 \supseteq C_0$ sit in $B \supseteq C$ with the unit-preserving way and that $\theta(C_0) \subseteq B_0$ for all $\theta \in \Theta$. We further assume that the restrictions of $E_C^B$ and $E_{\theta(C)}^B$ to $B_0$ give conditional expectations $E_{C_0}^{B_0}$ and $E_{\theta(C_0)}^{B_0}$ from $B_0$ onto $C_0$ and $\theta(C_0)$, respectively. Let us consider the following two reduced HNN extensions 
\begin{align*}
\left(A_0, E^{A_0}_{\pi\left(B_0\right)}\right) &= \left(B_0, E_{C_0}^{B_0}\right) \underset{C_0}{\bigstar} \left(\Theta|_{C_0}, \left\{E_{\theta(C_0)}^{B_0}\right\}_{\theta \in \Theta}\right), \\
\left(A, E^A_{\pi(B)}\right) &= \left(B, E_C^B\right) \underset{C}{\bigstar} \left(\Theta, \left\{E_{\theta(C)}^B\right\}_{\theta \in \Theta}\right) 
\end{align*}
with $\Theta|_{C_0} := \left\{ \theta|_{C_0}\ :\ \theta \in \Theta \right\}$. The reduced free products with amalgamations appeared in the procedure of construction are denoted by ${\mathcal B}_0$ and ${\mathcal B}$, respectively. In this setting, it is natural to ask when the following natural embedding exists:  
\begin{equation*}
A_0 \hookrightarrow A\ \text{by}\ 
\begin{cases}
\ b \in B_0 \mapsto b \in B; \\  
u\left(\theta|_{C_0}\right) \mapsto u(\theta),\ \theta \in \Theta. 
\end{cases}
\end{equation*}
To this question, we have a satisfactory answer as simple application of Blanchard and Dykema's work \cite{blanchard+dykema:PJM}. Namely, if all given conditional expectations have the faithful GNS representations, then there is such an embedding in the amalgamated free product level, ${\mathcal B}_0 \hookrightarrow {\mathcal B}$, and hence it is plain to see that $A_0$ is embedded into $A$ in the above-mentioned way. 

\subsection{Exactness of Reduced HNN Extensions}
Our construction has another advantage, which is a criterion for exactness.  To explain it, we should first remark that Dykema-Shlyakhtenko's result \cite[Proposition 4.1]{DykemaDima:procedinb2001} is still valid without any essential change even when the embedding maps $\iota_1$, $\iota_2$ are imposed upon the construction of free products with amalgamations. In fact,  when we consider the free product of unital $C^*$-algebras $A_1$, $A_2$ with amalgamation over a unital $C^*$-algebra $C$ via unital embedding maps $\iota_1 : C \rightarrow A_1$, $\iota_2 : C \rightarrow A_2$ with respect to conditional expectations $E_1 : A_1 \rightarrow \iota_1(C)$, $E_2 : A_2 \rightarrow \iota_2(C)$, it suffices only to replace, in their proof, the $C^*$-subalgebra $D$ of the $C^*$-algebra $A :=A_1 \oplus A_2$ and the completely positive map $\eta : A \rightarrow A$ by  
\begin{equation*} 
D :=\iota_1(C)\oplus\iota_2(C), \quad \eta\left(a_1, a_2\right) :=\left(\iota_1\circ\iota_2^{-1}\circ E_2\left(a_2\right),  \iota_2\circ\iota_1^{-1}\circ E_1\left(a_1\right)\right), 
\end{equation*}
respectively. 

Assume that a given unital $C^*$-algebra $B$ is exact and moreover that given conditional expectations $E_C^B : B \rightarrow C$, $E_{\theta(C)} : B \rightarrow \theta(C)$, $\theta \in \Theta$, have the faithful GNS representations. If $\Theta$ is a finite set, then Dykema and Shlyakhtenko's result \cite[Corollary 4.2]{DykemaDima:procedinb2001} implies that the reduced free product with amalgamation 
\begin{equation*}
\left({\mathcal B}, {\mathcal E}\right) = 
\left(B\otimes B\left(\ell^2\left(\Theta_1\right)\right), E_{\Theta} : \iota_{\Theta}\right) 
\underset{C\otimes\ell^{\infty}\left(\Theta_1\right)}{\bigstar} 
\left(B\otimes B\left(\ell^2\left(\Theta_1\right)\right), E_1 : \iota_1\right)     
\end{equation*} 
is exact since $B\otimes B\left(\ell^2\left(\Theta_1\right)\right)$ is clearly exact, and so is the reduced HNN extension 
\begin{equation*}
\left(A, E^A_{\pi(B)}\right) = \left(B, E_C^B\right) \underset{C}{\bigstar} 
\left(\Theta, \left\{E_{\theta(C)}^B\right\}_{\theta \in \Theta}
\right) 
\end{equation*}
too thanks to \cite[Proposition 7.1, (i)]{kirchberg1994}. When $\Theta$ is an infinite set, the reduced HNN extension $A = B\bigstar_C\Theta$ is still exact since $A=\varinjlim B\bigstar_C\Xi$ with finite subsets $\Xi \nearrow \Theta$ and \cite[Proposition 7.1, (iv)]{kirchberg1994}. On the other hand, if $A = B \underset{C}{\bigstar} \Theta$ is exact, then so should be $B$. Hence, the exactness of $B$ is necessary and sufficient for that of $A = B \underset{C}{\bigstar} \Theta$. With Remark 7.2, this fact in particular says that if a given countable discrete group $G$ is $C^*$-exact, then so is every HNN extension $G *_H \theta$ thanks to \cite[Theorem 5.2]{kirchbergwassermann1995}. This is indeed a fact mentioned in  \cite{guentner:procAMS2001}.


\begin{thebibliography}{99}

\bibitem{blanchard+dykema:PJM} 
E.~F.~Blanchard and K.~J.~Dykema, 
{\it Embeddings of reduced free products of operator algebras},  
Pacific J.~Math., {\bf 199} (2001), 1--19.  

\bibitem{connes:annals76} A.~Connes, 
{\it Classification of injective factors, Cases II$_1$, II$_{\infty}$, III$_{\lambda}$, $\lambda\neq1$}, 
Ann.~Math., {\bf 104} (1976), 73--115. 

\bibitem{connes:jot80} 
A.~Connes, 
{\it A factor of type ${\rm II}\sb{1}$ with countable fundamental group}, 
J.~Operator Theory, {\bf 4} (1980), 151--153.

\bibitem{connes-jones:bull-london} 
A.~Connes and V.~F.~R.~Jones, 
{\it Property $T$ for von Neumann algebras}, 
Bull.~London Math.~Soc., {\bf 17} (1985), 57--62. 

\bibitem{cowling-haagerup:invent89} 
M.~Cowling and U.~Haagerup, 
{\it Completely bounded multipliers and Fourier algebra of a simple Lie group of real rank one}, 
Invent.~math., {\bf 96} (1989), 507--549. 

\bibitem{DykemaDima:procedinb2001} 
K.~J.~Dykema and S.~Shlyakhtenko, 
{\it Exactness of Cuntz-Pimsner $C\sp *$-algebras}, 
Proc.~Edinb.~Math.~Soc., (2) {\bf 44} (2001), 425--444. 

\bibitem{falconetakesaki:quantumflow}
T.~Falcone and M. Takesaki, 
{\it The non-commutative flow of weights on a von Neumann algebra},  
J.~Funct.~Anal., {\bf 182} (2001), 170--206.

\bibitem{feldmanmoore:transAMS}
J.~Feldman and C.~C.~Moore, 
{\it Ergodic equivalence relations, cohomology, and von Neumann algebras, I, II}, 
Trans.~Amer.~Math.~Soc., {\bf 234} (1977), 289--324, 325--359.

\bibitem{gaboriau:invent2000}
D.~Gaboriau, 
{\it Co\^{u}t des relations d'\'{e}quivalence et des groupes},  
Invent.~Math., {\bf 139} (2000), 41--98.

\bibitem{guentner:procAMS2001}
E.~Guentner, 
{\it Exactness of the one relator groups}, 
Proc.~Amer.~Math.~Soc., {\bf 130} (2002), 1087--1093

\bibitem{haagerup:invent79}
U.~Haagerup, 
{\it An example of non-nuclear $C^*$-algebra which has the metric approximation property}, 
Invent.~math., {\bf 50} (1979), 279--293. 

\bibitem{hiaiueda:jot}
F.~Hiai and Y.~Ueda, 
{\it Automorphisms of free product type and their applications},
J.~Operator Theory, {\bf 50} (2003), 119--130.

\bibitem{higmanneumann^2:JLondonMath1949} 
G.~Higman, B.~H.~Neumann and H.~Neumann, 
{\it Embedding theorem for groups}, 
J.~London Math.~Soc., {\bf 24} (1949), 247--254.

\bibitem{kirchberg1994}
E.~Kirchberg, 
{\it Commutants of unitaries in UHF algebras and functorial properties of exactness}, 
J.~reine angew.~Math., {\bf 452} (1994), 39--77. 

\bibitem{kirchbergwassermann1995}
E.~Kirchberg and S.~Wassermann, 
{\it Operations on continuous bundles of $C^*$-algebras}, 
Math.~Ann., {\bf 303} (1995), 677--697. 

\bibitem{lyndonschupp:book}
R.C.~Lyndon and P.E.~Schupp, 
{\it Combinatorial Group Theory}, 
Classics in Mathematics. (Springer-Verlag, Berlin, 2001). xiv+339 pp. 

\bibitem{mcduff:uncountableexamples}
D.~McDuff, 
{\it Uncountably many ${\rm II}\sb{1}$ factors},  
Ann.~Math., (2) {\bf 90} (1969), 372--377.

\bibitem{murrayvonneumann4}
F.~Murray and J.~von Neumann, 
{\it On rings of operators, IV},
Ann.~Math., {\bf 44} (1943), 716--808. 

\bibitem{neshveyevstormer:jfa02}
S.~Neshveyev and E.~St{\o}rmer, 
{\it Ergodic theory and maximal abelian subalgebras of the hyperfinite factor}, 
J.~Funct.~Anal., {\bf 195} (2002), 239--261.

\bibitem{nielsen:jfa70}
O.~A.~Nielsen, 
{\it Maximal abelian subalgebras of hyperfinite factors. II}, 
J.~Funct.~Anal., {\bf 6} (1970), 192--202.

\bibitem{ocneanu:LNM}
A.~Ocneanu, 
{\it Actions of Discrete Amenable Groups on von Neumann Algebras}, 
Lecture Notes in Math., {\bf 1138}, 1985. 

\bibitem{ozawa:procAMS}
N.~Ozawa, 
{\it There is no separable universal II$_1$-factor}, 
Proc.~Amer.~Math.~Soc., {\bf 132} (2004), 487--490.

\bibitem{ozawa:acta}
N.~Ozawa, 
{\it Solid von Neumann algebras}, 
Acta Math., to appear. 

\bibitem{ozawapopa:invent}
N.~Ozawa and S.~Popa, 
{\it Some prime factorization results for type II$_1$ factors},
Invent.~Math., {\bf 156} (2004), 223--234. 

\bibitem{phillips:duke76}
J.~Phillips, 
{\it Automorphisms of full II$_1$ factors, with applications to factors of type III}, 
Duke Math. J., {\bf 43} (1976), 375--385. 

\bibitem{popa:jot83}
S.~Popa, 
{\it Orthogonal pairs of $*$-subalgebras in finite von Neumann algebras}, 
J.~Operator Theory, {\bf 9} (1983), 253--268. 

\bibitem{popa:adv}
S.~Popa, 
{\it Maximal injective subalgebras in factors associated with free groups}, 
Adv.~Math., {\bf 50} (1983), 27--48.

\bibitem{popa:invent:irredindex>4}
S.~Popa, 
{\it Markov traces on universal Jones algebras and subfactors of finite index}, 
Invent.~Math., {\bf 111} (1993), 375--405. 

\bibitem{popa:betti}
S.~Popa, 
{\it On a class of type II$_1$ factors with Betti numbers invariants}, 
Preprint (2002).  

\bibitem{popashlyakhtenko:preprint}
S.~Popa and D.~Shlyakhtenko, 
{\it Universal properties of $L\left({\mathbb F}_{\infty}\right)$ in subfactor theory}, 
Acta Math., {\bf 191} (2003), 225--257.  

\bibitem{radulescu:jams}
F.~R\u{a}dulescu, 
{\it The fundamental group of the von Neumann algebra of a free group with infinitely many generators is ${\mathbf R}_+\setminus\{0\}$},  
J.~Amer.~Math.~Soc., {\bf 5} (1992), 517--532.

\bibitem{radulescu:freegroupsubfactor}
F.~R\u{a}dulescu, 
{\it Random matrices, amalgamated free products and subfactors of the von Neumann algebra of a free group, of noninteger index}, 
Invent.~Math., {\bf 115} (1994), 347--389. 

\bibitem{sekine:mathscand}
Y.~Sekine, 
{\it Cartan subalgebras in fixed point algebras of finite group actions}, 
Math. Scand., {\bf 70} (1992), 281--292. 

\bibitem{serre:book}
J.~P.~Serre, 
{\it Trees} (Translated from the French by John Stillwell), Springer-Verlag, Berlin-New York, 1980. 

\bibitem{shlyakhtenkoueda:crelle}
D.~Shlyakhtenko and Y.~Ueda, 
{\it Irreducible subfactors of $L\left({\mathbb F}_{\infty}\right)$ of index $\lambda > 4$}, 
J.~reine angew.~Math., {\bf 548} (2002), 149 -- 166.  

\bibitem{takesaki:book}
M.~Takesaki, 
{\it Theory of Operator Algebras, I, II, III}, 
Encyclopaedia of Mathematical Sciences. Vol.~{\bf 124} (2002), Vol.~{\bf 125} (2003), Vol.~{\bf 127} (2003),  Operator Algebras and Non-commutative Geometry, {\bf 5}, {\bf 6}, {\bf 8}, Springer-Verlag, Berlin.

\bibitem{ueda:journalofmathsocJpn}
Y.~Ueda, 
{\it A minimal action of the compact quantum group $\mathrm{SU}_q(n)$ on a full factor}, 
J.~Math.~Soc.~Japan, {\bf 51} (1999), 449--461. 

\bibitem{ueda:pacific}
Y.~Ueda, 
{\it Amalgamated free product over Cartan subalgebra}, 
Pacific J.~Math., {\bf 191} (1999), No.2, 359--392. 

\bibitem{ueda:us-jpn}
Y.~Ueda, 
{\it Amalgamated free product over Cartan subalgebra, II: Supplementary results \& examples}, 
Advanced Studies in Pure Mathematics, {\bf 38} (2004) ``Operator Algebras and Applications" 239-265. 

\bibitem{ueda:transAMS}
Y.~Ueda, 
{\it Fullness, Connes' $\chi$-groups, and ultra-products of 
amalgamated free products over Cartan subalgebras}, 
Trans.~Amer.~Math.~Soc., {\bf 355} No.1 (2003), 349--371. 

\bibitem{voiculescu:lecturenotes1132}
D.~Voiculescu, 
{\it Symmetries of some reduced free product $C^*$-algebras}, 
Lecture Notes in Math., {\bf 1132} (1984), 556--588.  

\bibitem{voiculescu:progress}
D.~Voiculescu, 
{\it Circular and semicircular systems and free product factors}, 
in Operator Algebras, Unitary Representations, Enveloping Algebras, and Invariant Theory (Paris, 1989), 45--60, Progr. Math., {\bf 92}, Birkh\"{u}ser Boston, Boston, MA, 1990. 

\bibitem{voiculescu:bulllondon2002}
D.~Voiculescu, 
{\it Free entropy}, 
Bull.~London Math.~Soc., {\bf 34} (2002), 257--278.

\bibitem{voiculescudykemanica:monograph}
D.~Voiculescu, K.~Dykema and A.~Nica, 
{\it Free Random Variables},  
CRM Monograph Series, {\bf 1}. (Amer. Math. Soc., Providence, RI, 1992).

\bibitem{yamagami:aspect-modular}
S.~Yamagami, 
{\it Algebraic aspects in modular theory}, 
Publ.~Res.~Inst.~Math.~Sci., {\bf 28} (1992), 1075--1106.

\end{thebibliography}
\end{document}